\documentclass[a4paper,twoside,english]{amsart}
\usepackage[T1]{fontenc}
\usepackage[latin9]{inputenc}
\usepackage{babel}
\usepackage{verbatim}
\usepackage{amsthm}
\usepackage{amsmath}
\usepackage{amssymb}
\usepackage[unicode=true,
 bookmarks=true,bookmarksnumbered=true,bookmarksopen=false,
 breaklinks=false,pdfborder={0 0 1},backref=false,colorlinks=false]
 {hyperref}
\hypersetup{
 pdfauthor={D. Masser  and U. Zannier}}

\makeatletter

\pdfpageheight\paperheight 
\pdfpagewidth\paperwidth

\theoremstyle{plain}

 \newcommand\thmsname{\protect\theoremname}
 \newcommand\nm@thmtype{theorem}
 \theoremstyle{plain}

  \theoremstyle{remark}
  \newtheorem{rem}[thm]{\protect\remarkname}
  \theoremstyle{definition}
  \newtheorem*{example*}{\protect\examplename}
  \theoremstyle{definition}
  
  \theoremstyle{plain}
  
  \theoremstyle{plain}
  
  \theoremstyle{plain}

\makeatletter
  \theoremstyle{definition}
  \newtheorem{my@rem}[thm]{Remark}
  
\makeatother

\makeatother

  \providecommand{\examplename}{Example}
  \providecommand{\lemmaname}{Lemma}
  \providecommand{\propositionname}{Proposition}
  \providecommand{\remarkname}{Remark}
  \providecommand{\theoremname}{Theorem}
\providecommand{\theoremname}{Theorem}
 \providecommand{\corollaryname}{Corollary}

\usepackage[left=3.1cm,right=3.5cm,top=3cm,bottom=3cm]{geometry}

\usepackage{mathtools}

\def\CVD{{\hfill\hfil{\lower 2pt\hbox{\vrule\vbox to 7pt
{\hrule width  5pt\varphifill\hrule}\varphirule}}}\par}

\begin{document}
\title{ALAN BAKER}

\author{David Masser}
\address{Department of Mathematics and Computer Science\newline
\indent University of Basel\newline
\indent Spiegelgasse 1, 4051 Basel, Switzerland}
\email{david.masser@unibas.ch}

\begin{abstract}
Alan Baker, Fields Medallist, died on 4th February 2018 in Cambridge England after a severe stroke a few days earlier. 

In 1970 he was awarded the Fields Medal at the International Congress in Nice on the basis of his outstanding work on linear forms in logarithms and its consequences. Since then he received many honours including the prestigious Adams Prize of Cambridge University, the election to the Royal Society (1973) and the Academia Europaea; and he was made an honorary fellow of University College London, a foreign fellow of the Indian Academy of Science, a foreign fellow of the National Academy of Sciences India, an honorary member of the Hungarian Academy of Sciences, and a fellow of the American Mathematical Society.

In this article we survey Alan Baker's achievements.

\end{abstract}

\maketitle

\section{Life and Career}\label{sec1}
Alan's paternal grandparents were known as Marks and Mathilda Backer, married 1902 in Lithuania. The surname  presumably changed when children arrived. His parents were Barnet and Bessie (n\'ee Sohn). Alan was born on 19th August 1939 into this Jewish family. His earliest memories of wartime England were of evacuation to Camberley, Surrey. After the war the family moved to Forest Gate in East London, where he spent most of his early life. From a very early age he showed signs of mathematical brilliance (see the comments later about ``brainbox'') and was encouraged by his parents. Already his father (who had been at school with Jacob Bronowksi) was very gifted in this direction, but did not have the opportunity to develop, and became a tailor instead. That may have explained Alan's clothes sense; he was always well turned out with quality suits and tasteful ties (which were however not always entirely appropriate to his later travelling, for example on the beach at Nice after being awarded the Fields Medal or climbing a hill in the Australian Bush among snakes - such episodes may have contributed to his later investing in a distinctive yellow sun hat).

Alan's first education was at a Franciscan Convent (London E.7) 1945-49, followed by Godwin County Primary School 1949-50 and Stratford Grammar School 1950-58. He then went with a State Scholarship to University College London (he recalls that the staff included J.W. Archbold, L.S. Bosanquet, P. Du Val, T. Estermann, C.A. Rogers, and K. Roth), where he studied mathematics 1958-61, obtaining a first class honours B. Sc. (Special) degree.  

He then moved to Trinity College (where he would soon be based for the rest of his life) in Cambridge to study 1961-64 with Harold Davenport, one of the leading number theorists at the time with many international connections, of whom Alan writes: 

\medskip
``{\it An excellent mathematician from whom I learnt a great deal; but I tended to follow my own lines of research}''. 

\medskip
\noindent
Still, Chris Morley recalls a loud mathematical conversation between them in the Trinity Parlour, and they did write a paper together (see later). It seems that in later life he enjoyed imitating Davenport's Lancashire (Accrington) accent (and the writer has carried on this worthy tradition of imitating one's supervisor). He obtained his Ph.D. 1965 and M.A. 1966, by which time he had already been awarded a Prize Fellowship 1964-68 also at Trinity. During this period, which also included a year 1964-65 at University College, he took on John Coates as a Ph.D. student and recalls also sharing responsibilities with Davenport in supervising T.W. Cusick, M.N. Huxley, H.L. Montgomery, and R.W.K. Odoni. 

In 1966 he was appointed as Assistant Lecturer, in 1968 as Lecturer, in 1972 promoted to Reader in the Theory of Numbers, and finally in 1974 he was elected to a personal chair for Pure Mathematics, all at Cambridge. Apparently he liked to draw attention to the unusual chronological order of Fields Medal then Fellow of Royal Society then Professor. During that period he supervised for Ph.D. also, in chronological order, the writer, Cameron Stewart, Yuval Flicker, Roger Heath-Brown, Richard Mason, Mark Coleman and E.J. Lee (whose full name I have been unable to find out). At the time of writing, the Mathematical Genealogy Project lists 463 descendants. 

The writer recalls meeting Alan every two weeks or so, in his college rooms in Whewell's Court opposite the Great Gate. Actually these rooms (fitting C.P. Snow's description ``not specially agreeable'') had been occupied previously by G.H. Hardy (about which the great German mathematician David Hilbert - see later - indignantly wrote to the Master that Hardy was the best mathematician, not only in Trinity, but in England, and should therefore have the best rooms). One had to ascend a slightly low and narrow spiral staircase, and then knock on the thick wooden door. Often in College {\it une porte peut en cacher une autre} and probably he did indeed have a second door immediately behind; at any rate he was sometimes a long time coming and I supposed that he did not hear me through the two doors (also against the heavy traffic then in Trinity Street). On each visit I felt obliged to knock or bang harder and yet harder, sometimes bringing with me a heavy book to protect my hands. 

As mentioned, Alan Baker was firmly based in Cambridge; it seems that college life there suited him especially in the style of Trinity, whose society he enriched for many years in an unspectacular way, for example being an interesting conversationalist (despite the impression given by some newspaper obituaries). He also threw after-seminar parties fuelled by particularly strong beer (on one occasion Frank Adams challenged others to climb around a table, underneath it, without the feet touching the ground, after he had himself demonstrated that it could be done). In the Combination Room there was a bottle of Madeira kept specially for him. He was a reasonable ballroom dancer, enthusiastically participating in the College May Ball. He played regularly on the Trinity Bowling Green, using an unconventional throw which delivered the ball from waist level instead of lower down. 

Outside the College, he enjoyed playing the slot machines in pubs, or playing table-tennis and snooker at the Graduate Centre. Eva McLean (n\'ee Gordon) writes of one occasion:

\medskip
``{\it In 1975 a vicious rapist was terrorising the female population of Cambridge. One evening when I was leaving the Graduate Centre one of the porters who knew me well expressed great concern that I would be walking home alone. Not to worry, I assured him brightly, pointing to Alan who was going to escort me. The burly ex-policeman looked him up and down, all five foot six of him, grabbed his coat and joined us on our way.}''

\medskip
\noindent
and draws attention to existential dangers of a different kind:

\medskip
``{\it Alan was once barricaded in his flat for several days and would not answer the door as a collector, all the way from America, kept coming back pleading, in vain, for a contribution to his sperm bank.}''

\medskip
He had learnt to drive in America, and in Cambridge he bought a Rover car, apparently as an investment (``{\it you understand}'', as he insisted), although Eva McLean soon explained to him that there were much better ways of investing (and he seems to have acted very well on this advice). He drove with enthusiasm despite colliding with a stationary fire engine on his maiden voyage, and despite being stopped for speeding which caused him to arrive late in chairing a session in London. William Chen recalls  that it was actually the inaugural lecture of a colleague, and that Alan was very eager to explain that he had not exceeded 90mph, so he was only speeding rather than driving recklessly. On one occasion the car got badly stuck in a driveway; Alan however remained cool and concentrated and succeeded in extricating it without damage. Later on the enthusiasm waned and a rusting hulk with flat tyres had to be removed from New Court. He had a flat in Hendon and enjoyed life there in London too, for example upmarket restaurants, or the theatre at which he always bought the best seats (as at the Footlights Theatre in Cambridge). The latter interest may be due to his cousin Heather Rechtman, who writes: 

\medskip
``{\it We lived near Stratford East (London) where the fledgeling company Theatre Workshop was just beginning to make a name for itself and as a stage struck teenager I dragged him along to many of their plays}''.  

\medskip
\noindent
That was of course edgy stuff under Joan Littlewood; but he also enjoyed musicals. About one such trip to the theatre Eva McLean narrates: 

\medskip``{\it ... we stopped over at his Hendon flat. While there he decided to show me his mother's mink coat. When he opened the wardrobe door hundreds of moths awoke to greet us. His response to this crisis? He shut the door firmly.}''

\medskip
\noindent
as well as:

\medskip ``{\it On another eventful outing to the West End, we arrived in good time and right there by the theatre came across a perfect parking place. Alan, however, decided it was too good to be true and so the search continued. When we next passed the spot, the parking space was gone - as was the first act of the play by the time we were inside.}''

\medskip
Her funeral tribute (23rd February 2018) sums up such anecdotes:

\medskip
``{\it So, in many ways, he was what Americans term `just a regular guy'.}''

\medskip
In his private life he was quite relaxed, never going beyond sometimes stamping his foot when frustrated. However, in his professional life he was often reserved, even insecure in some ways, and could on occasion be difficult. There were episodes, not just in Cambridge, which perhaps still have not been entirely forgotten. But it is now impossible to give any balanced accounts.

\medskip
He was enthusiastic about travel, which in America already started in 1969 with visiting professorships in Ann Arbor Michigan and Boulder Colorado and in 1970 membership of the Institute for Advanced Study in Princeton New Jersey. He had at least three offers of chairs. One of these he turned down because the ivy-clad walls were too much like those in Cambridge, as was the weather. He preferred both modernity and sunshine. It was also thought that he wanted to look after his mother in London (this writer recalls around 1980 seeing two people slowly crossing a Maryland street, and recognizing one but not the other, and he followed them into a restaurant where Alan introduced his mother). He made many American visits, and Eva McLean has shown me a lot of his postcards and other correspondence from abroad, for example from Texas (1984), New York (1987), and California (1989). Of course as his reputation grew he was able to widen his horizons, for example to Australia, China/Hong Kong (missive 1988), India, Japan (postcard 1983), Russia, and many parts of Europe.

In later life he made regular trips to Switzerland (postcards 1989 and 1990) to work with Gisbert W\"ustholz at ETH (a 1988 missive opens with ``greetings from a gnome of Z\"urich''). There were social events too; and he thought nothing of occasionally bringing presents such as a frying pan to dinner invitations - unconventional to be sure, but with its own unassailable logic. It was there, during a conference in honour of his 60th birthday, that he gave an entertaining and surprisingly candid speech about his life, starting with his recollections of wartime London, also mentioning that he was regarded as the ``{\it brainbox}'' of the family, including one more rendering of Davenport's accent; and ending with his regrets about never marrying. Some people already knew that these regrets were not abstract but concerned specific ladies. I quote again from the funeral tribute:

\medskip
``{\it I first met Alan in 1975 playing table tennis in the Graduate Centre. I was immediately taken by his generosity of spirit. Unlike some of the others, he never minded being beaten. He proved to be equally generous and attentive as a suitor, and later as an old friend.}''

\medskip
His last years were made more difficult by increasing deafness and a series of falls (he did not pay serious attention to medical advice on these). He had long since advanced from Whewell's Court to Great Court and indeed had been proud to have the very best set of rooms there, like a maisonette on two floors. It overlooked the Bowling Green where he now could no longer play - see the cover of \cite{BG}, the proceedings of the Z\"urich conference.

During a Feast at Trinity in his memory he was described as ``idiosyncratic'' in place of words like ``eccentric'' or ``enigmatic''.

\medskip
Other biographical articles about Alan Baker can be found in the Hardy-Ramanujan Journal \cite{HR}.

\medskip
See also \cite{AA} for a scientific appraisal by Gisbert W\"ustholz, which also contains a complete list of his publications.

\section{Mathematics - a preview}\label{sec2}

For a more structured narrative we shall divide Alan Baker's work into eight categories (and in scientific detachment usually drop the first name):

\bigskip
(a) Diophantine approximation.

\medskip
(b) Linear forms in logarithms.

\medskip
(c) Diophantine equations.

\medskip
(d) Elliptic functions.

\medskip
(e) Class numbers.

\medskip
(f) Abcology.

\medskip
(g) Miscellaneous.

\medskip
(h) Books.

\bigskip
But before starting, we would like to mention perhaps the most easily stated of all his deep results. 

\medskip
When we make a list of the perfect squares
$$1,4,9,16,25,36,49,64,81,100,121,144,\ldots,143384152921,\ldots$$
we see that the gaps between consecutive members get larger and larger (and in a regular way). Similarly for the perfect cubes
$$1,8,27,64,125,216,343,512,729,1000,1331,1728,\ldots,143384152904,\ldots.$$
But if we mesh the two lists together to get the ``squbes''
$$1,4,8,9,16,25,27,36,49,64,81,\ldots,143384152904,143384152921,\ldots$$
then it is not so clear that the gaps get large. Indeed it was proved by L.J. Mordell only in 1922; thus for example the gap 17 (twice above) or the gap 1621 (to take a year apparently at random) occurs at most finitely often. Unfortunately Mordell's proof gave no way of determining all the occurrences of a given gap.

To do this amounts to specifying a non-zero $k$ in the set $\bf Z$ of rational integers and finding all $x,y$ in $\bf Z$ with
\begin{equation}\label{mord1}
y^2=x^3+k.
\end{equation}
Baker achieved this in 1968 by showing that they all satisfy
\begin{equation}\label{mord2}
\max\{|x|,|y|\} ~\leq~ \exp\{(10^{10}|k|)^{10000}\}.
\end{equation}
Despite (\ref{mord1}) being around since at least the year 1621 (and the cases $k=-2,-4$ were set by Fermat in 1657 as a challenge to ``you English''), there were no estimates at all for $x,y$ until (\ref{mord2}) nearly 350 years later. 

Thus for example to find all gaps 1621 one just has to examine all $y^2,x^3$ with $y,x$ between 1 and say
$$10^{10^{132098}}.$$
This looks hopelessly impractical; yet we will see later that Baker (with Davenport) found exceedingly efficient ways to do such things. 

There is an attractive single-sentence reformulation: for all positive integers $x,y$ with $x^3 \neq y^2$ we have
$$|x^3-y^2| > 10^{-10}(\log x)^{1/10000}.$$

\section{Diophantine approximation}\label{sec3}

It is classical that
$$\pi=3.1415926\ldots$$
and 
$${355 \over 113}=3.1415929\ldots$$
are suspiciously close. This raises the natural question: given a real number $\xi$, how well can we approximate it by a rational number $p/q$? An answer in convenient form was found by Dirichlet: provided $\xi$ is not already rational, we can find infinitely many $p/q$ with
\begin{equation}\label{ub}
\left|\xi - {p \over q}\right| < {1 \over q^2}.
\end{equation} 

We pause to give the simple proof, which involves the Box Principle or Pigeonhole Principle. Pick any $Q \geq1$ in $\bf Z$. For $i=0,1,\ldots,Q$ we can find $p_i$ in $\bf Z$ with the $Q+1$ pigeons $\theta_i=i\xi-p_i$ in the interval from 0 to 1. We divide this interval into $Q$ holes of length $1/Q$. There are more pigeons than holes, so at least one hole must contain at least two pigeons. With say $\theta_j$ and $\theta_k$ ($j<k$) this leads at once to
\begin{equation}\label{dir}
|q\xi-p| \leq {1 \over Q}
\end{equation}
for $p=p_k-p_j$ and $q=k-j$ satisfying $1 \leq q \leq Q$. And now (\ref{ub}) follows; it is not too hard to see that we get infinitely many $p/q$ as $Q$ varies.

\medskip
But the answer to the next natural question of whether we can beat (\ref{ub}) depends critically on the number $\xi$.

For a class of $\xi$ particularly interesting to number-theorists, this topic can be said to have begun in earnest with Liouville in 1844, although it underlies the older concept of continued fractions such as
$$3+{1\over {7+{1 \over 16}}}={355\over113}$$
or infinite ones like
$$3+{1\over {7+{1 \over 15+{1 \over 1+{1 \over 292+\cdots}}}}}=\pi.$$
More generally
\begin{equation}\label{cf1}
a_0+{1\over {a_1+{1 \over a_2+{1 \over a_3+{1 \over a_4+\cdots}}}}}
\end{equation}
is usually linearized to just
\begin{equation}\label{cf3}
x=[a_0;a_1,a_2,a_3,a_4,\ldots].
\end{equation}
For example
\begin{equation}\label{cf4}
\alpha=[1;1,1,1,1,\ldots]={1+\sqrt{5} \over 2},
\end{equation}
and this is actually an algebraic number in that $\alpha^2-\alpha-1=0$. On the other hand Liouville's work (see later) shows that if the positive integers $a_0,a_1,a_2,a_3,a_4,\ldots$ increase very rapidly, then $x$ in (\ref{cf3}) does not satisfy any equation 
\begin{equation}\label{alg5}
b_0x^d+b_1x^{d-1}+\cdots+b_d=0
\end{equation}
for $b_0,b_1,\ldots,b_d$ in $\bf Z$ not all zero. Thus by definition $x$ is a transcendental number.

In 1906 Maillet had given a different sort of transcendental continued fraction. A typical example is to take (\ref{cf4}) and replace 1 by 2 in the positions $k_1,k_2,\ldots$, where now these $k_1,k_2,\ldots$ increase very rapidly. 

\medskip
In his very first paper \cite{1w} from 1962 Baker simplified and improved that work and also made the estimates more explicit. A consequence here is that it suffices to take $k_n$ as small as $4^n$ for transcendence. For the proof one notes that something like
$$x=[1;1,1,1,1,1,1,1,1,1,1,2,1,\ldots]$$ 
is rather close to $\alpha$ above, with similar approximations further along. In this case an easy generalization of Liouville is applicable, but for other examples it is necessary to use Klaus Roth's ``revolutionary improvement'' (Ian Cassels) of Liouville's result, or more precisely a consequence by Davenport and Roth, and even LeVeque's generalization of Roth.

Such acknowledgements as Baker's 

\medskip
``{\it I should like to thank Professor Davenport for his valuable suggestions and help in preparing the manuscript.}''

\medskip 
\noindent
can be seen quite often in the journal Mathematika around this time (which Davenport founded).

\medskip
The paper \cite{4w} from 1964 can be considered as a sort of continuation of \cite{1w}. To describe some of its results we must recall that what Liouville proved is that for any algebraic number $\alpha$ of degree $d \geq 2$ (the smallest integer such that $\alpha$ satisfies an equation (\ref{alg5}) above) there is $c>0$ such that
\begin{equation}\label{liou6}
\left|\alpha - {p \over q}\right| \geq {c \over q^d}
\end{equation}
for all $p$ and $q \geq 1$ in $\bf Z$. When $d \geq 3$ Roth improved this to
\begin{equation}\label{roth7}
\left|\alpha - {p \over q}\right| \geq {c \over q^\kappa}
\end{equation}
for any $\kappa>2$, where now $c$ is allowed to depend on $\kappa$. The Box Principle as in (\ref{ub}) shows that this is essentially best possible.

Equivalently if for some $\xi$ and $\kappa>2$ there is an infinite sequence of $p_i$ and $q_i \geq 1$ in $\bf Z$ with $p_i/q_i$ different and 
\begin{equation}\label{app8}
\left|\xi - {p_i \over q_i}\right| < {1 \over q_i^\kappa}~~~~~~~(i=1,2,\ldots)
\end{equation} 
then $\xi$ must be transcendental. Baker then shows under an additional condition that $\xi$ cannot be too close to an algebraic number in the sense of what is since called a ``transcendence measure''. The condition is that there should exist $\lambda$ with $q_{i+1} \leq q_i^\lambda$ for all $i$. In that case he shows that for each $n$ there is $\mu_n$ (possibly depending on $\xi$) such that
\begin{equation}\label{tm9}
|h_0\xi^n+h_1\xi^{n-1}+\cdots+h_n| > H^{-\mu_n}
\end{equation}
for all $h_0,h_1,\ldots,h_n$ in $\bf Z$ not all zero, where $H=\max\{2,|h_0|,|h_1|,\ldots,|h_n|\}$.

The significance of this is the following. From (\ref{tm9}) by definition (due to Kurt Mahler) $\xi$ cannot be a so-called $U$-number, even though (\ref{app8}) implies that the ``partial quotients'' $a_0,a_1,\ldots$ in (\ref{cf3}) for $x=\xi$ are unbounded. On the other hand the partial quotients for things like $\sum_{k=1}^\infty 2^{-k!}$ are also unbounded and it is a $U$-number (for example (\ref{tm9}) fails already for $n=1$ and $h_0=2^{k!}$ with $k$ large). Furthermore he shows that the method of \cite{1w} produces $U$-numbers with bounded partial quotients, and that a suitable generalization of the above result to quadratic fields produces $\xi$ that are not $U$-numbers but still with bounded partial quotients.

More succinctly, there is no correlation between the properties of having bounded partial quotients and being $U$-numbers. In fact Baker is able to sharpen this to involve $T$-numbers and $S$-numbers (see just below) also in Mahler's classification, and that was his main motivation. A spin-off is that either $T$-numbers exist or $S$-numbers ``of type exceeding 1'' exist. This foreshadows Wolfgang Schmidt's breakthrough four years later that $T$-numbers exist.

The proofs are rather formidable; indeed (\ref{app8}) and (\ref{tm9}) lie close to a strengthening of Roth's Theorem and accordingly Baker has to ramp up Roth's entire machinery, even in the situation of LeVeque's generalization.

This time he writes only

\medskip
``{\it I am indebted to Professor Davenport for valuable suggestions in connection with the present work.}''

\medskip
In \cite{10w} from 1967 he returned to these themes, with generalizations to several numbers at least in the case of bounded partial quotients. This property for a single $\xi$ is equivalent to $|\xi-p/q| \geq c/q^2$ analogous to (\ref{liou6}), and the natural extension to a pair $(\xi_1,\xi_2)$ is 
$$\max\left\{\left|\xi_1-{p_1 \over q}\right|,\left|\xi_2-{p_2 \over q}\right|\right\} \geq {c \over q^{3/2}}$$
for all $p_1,p_2$ and $q \geq 1$ in $\bf Z$. He shows for example that there are $U$-numbers $\xi_1,\xi_2$ satisfying this, and also $\xi_1,\xi_2$ that are not $U$-numbers. This time the proofs need no Roth-type considerations but follow a still intricate ``interval nesting'' technique of Cassels and Davenport.

\medskip
In 1966 Baker \cite{T} (published in the Proceedings of this Society) investigated certain ``metrical'' properties of $S$-numbers. These are defined by refining (\ref{tm9}) as follows. There is $\omega$ such that for each $n$ and each $\kappa>n\omega$ we can find $c>0$ (possibly depending on $\xi,n,\kappa$) such that
\begin{equation}\label{tm10}
|h_0\xi^n+h_1\xi^{n-1}+\cdots+h_n| > cH^{-\kappa}
\end{equation}
as in (\ref{tm9}). For simplicity we stick to real $\xi$. Much as before, the Box Principle shows that $\omega \geq 1$. It had been conjectured by Mahler in 1932 that in fact we can take $\omega=1$ (``of type exactly 1'') for almost all $\xi$ in the sense of Lebesgue measure. This was proved by Sprindzhuk in 1965. Having seen only preliminary announcements of that result, Baker was able to refine (\ref{tm10}) even further to things like
\begin{equation}\label{tm10a}
|h_0\xi^n+h_1\xi^{n-1}+\cdots+h_n| > cH^{-n}(\log H)^{-\lambda}
\end{equation}
for any $\lambda>n$.

\medskip
Later on in 1970 \cite{60w} with Schmidt he considered further refinements in terms of Hausdorff dimension; these are closer to Koksma's classification into $S^*$-,$T^*$- and $U^*$-numbers based instead on the distance from $\xi$ to algebraic numbers.

\medskip
Many of the results above can be expressed more concisely through the notation 
$$\|x\|=\min_{m \in {\bf Z}}|x-m|$$ 
for $x$ in the field $\bf R$ of real numbers, meaning the distance to the nearest integer. For example (\ref{liou6}) says that $\|q\alpha\| \geq c/q^{d-1}$ for all $q \geq 1$.

\medskip
This notation can be used in other contexts. Returning to 1964, we may cite his nice note \cite{3w}, which shows that
\begin{equation}\label{prod11}|Q|\|Q\Theta_1\|\|Q\Theta_2\| \geq e^{-5}
\end{equation}
for any non-zero polynomial $Q$ in say ${\bf R}[t]$, where $\Theta_1=e^{1/t},\Theta_2=e^{2/t}$ (same $e$ but that is irrelevant) are interpreted as formal power series in the topological completion with respect to the valuation $|t|=e$, and 
$$\|X\|=\min_{M \in {\bf R}[t]}|X-M|$$
measures the distance to the nearest polynomial. The significance here is that the analogue of (\ref{prod11}) with ${\bf R}[t]$ replaced by $\bf Z$ (see just below) is thought to be false, by a famous conjecture of J.E. Littlewood (still unsolved). The note makes explicit an earlier result of Davenport and Lewis and provides a simpler proof, essentially by differentiating and using a Pad\'e approximation (see later in this section). In fact the Pad\'e element can be eliminated to give a yet simpler proof as follows.

\medskip
It can be checked that 
$$\Delta=\left|\begin{array}{ccc}Q&P_1&P_2 \\ t^2Q'&t^2P_1'+P_1&t^2P_2'+2P_2 \\t^4Q''+2t^3Q'&t^4P_1''+(2t^3+2t^2)P_1'+P_1&t^4P_2''+(2t^3+4t^2)P_2'+4P_2\end{array}\right|$$
remains unchanged on replacing $P_1,P_2$ by ${\mathcal E}_1=P_1-Qe^{1/t},{\mathcal E}_2=P_2-Qe^{2/t}$ respectively.

Now suppose $Q,P_1,P_2$, even in ${\bf C}[t]$, are all non-zero. Inspecting the coefficients of smallest powers of $t$ in each entry shows that $\Delta \neq 0$. Thus $|\Delta| \geq 1$.

Next for $Q \neq 0$ choose $P_1,P_2$ so that $|{\mathcal E}_1|=\|Qe^{1/t}\|<1, |{\mathcal E}_2|=\|Qe^{2/t}\|<1$. Clearly $P_1 \neq 0, P_2 \neq 0$. Using $|X'| \leq e^{-1}|X|$ we get easily
$$1 \leq |\Delta| \leq e^3|Q||{\mathcal E}_1||{\mathcal E}_2|$$
giving a slight improvement of the result.

\medskip
He evidently had this note in mind when writing \cite{7w} in 1965, back in $\bf R$. A special result there is that for any $\epsilon>0$ there is $c>0$ such that
\begin{equation}\label{prod12}
q\|q\theta_1\|\|q\theta_2\| \geq {c\over q^\epsilon}
\end{equation}
for any $q \geq 1$ in $\bf Z$, where $\theta_1=e,\theta_2=e^2$.

More generally he treats a product of several terms (as indeed in \cite{3w} above) involving various numbers $\theta=e^\phi$ with (different) rational $\phi$ (that then implies quite easily the transcendence of $e$ with quite a good measure); and he is even able to refine the $q^\epsilon$. Of course one can no longer do things just by differentiation; but that does play a role and the much more elaborate proof uses techniques introduced by Carl Ludwig Siegel in his work on the so-called $E$-functions, together with so-called transference principles.

\medskip
The paper \cite{13w} from 1967 is Baker's only excursion into the Hardy-Littlewood Circle Method. Now one wants to know that certain approximations actually exist, and this time with prime values of the variables. For example, he shows that given any $\kappa$ there are infinitely many primes $p_1,p_2,p_3$ with
$$|p_1-\sqrt{2}p_2-\sqrt{3}p_3| < (\log P)^{-\kappa},$$
where $P=\max\{p_1,p_2,p_3\}$ (with extension to arbitrary coefficients). This was the first such result with an explicit function of $P$ on the right-hand side. As Bob Vaughan (who later greatly improved the result, to Baker's surprise and dismay) has pointed out, he was also aware that a ``localized'' version in the style of (\ref{dir}), with $P$ given in advance and just $P \geq \max\{p_1,p_2,p_3\}$, is in principle impossible. We may remark that the result in the enigmatic footnote on the second page of this \cite{13w} had previously been proved by Vinogradov.

\medskip
Returning again to 1964 we examine \cite{5w} in some detail. The starting point is again (\ref{liou6}) and (\ref{roth7}). There is a fundamental difference between these two results. Namely, given $\alpha$ in (\ref{liou6}), one can calculate an explicit value of $c$ (and in fact rather easily). One says informally that (\ref{liou6}) is ``effective''. 

In fact not so much attention had been paid to this concept, which dramatically increased in importance thanks to the works of Alan Baker. 

\medskip
On the other hand, given $\alpha$ and $\kappa$ in (\ref{roth7}), no-one knows how to calculate any value of $c$. This is true even for the simple-looking 
\begin{equation}\label{lb13}
\left|2^{1/3}  -  {p \over q}\right|  \geq {c \over q^{2.1}}.
\end{equation}
The obstacle is a logical twist already in Thue's improvement on (\ref{liou6}) in 1909. 

Namely, his argument works very well if there happens to be already a very good rational approximation $p_?/q_?$ to $\alpha$. This sounds paradoxical, but it has the effect of ``repelling'' other $p/q \neq p_?/q_?$ through the triangle inequality
\begin{equation}\label{tri}
\left|\alpha - {p \over q}\right| \geq \left|{p_? \over q_?} - {p \over q}\right|-\left|\alpha - {p_? \over q_?}\right|\geq {1 \over q_?q}-\left|\alpha - {p_? \over q_?}\right|.
\end{equation}
And if there are no such good approximations $p_?/q_?$ that is also good news. However we have no way of resolving this dichotomy; and in the first case $c$ will depend on the uncontrollable $p_?/q_?$.

In \cite{5w} Baker treats some numbers of this special form $\alpha= (a/b)^{1/d}$ for positive $a,b$ in $\bf Z$ as in (\ref{lb13}). He proves (\ref{roth7}), but with different quantifiers. Namely, for each such $\alpha$ he finds an explicit $\kappa$ such that (\ref{roth7}) holds for some effective $c>0$. Of course this already holds for $\kappa=d$ by (\ref{liou6}); but the new $\kappa$ can sometimes be less than $d$. Furthermore, and this is Baker's key observation, one can get $\kappa$ arbitrarily close to 2 by choosing $\alpha$ carefully, usually with $a,b$ approximately equal.  This amounts to imposing $p_?/q_?=1$ in advance. For example if $b=3^{316}$ then
\begin{equation}\label{lb14}
\left|\left({b+1 \over b}\right)^{1/3}  -  {p \over q}\right|  \geq {10^{-1024522} \over q^{2.1}},
\end{equation}
which is effective ``but only just'' (as Estermann's Scotsman \cite{est} said on p.33); and not as elegant as (\ref{lb13}).  

Since Thue's work it is known that any $\kappa<d$ leads to consequences for diophantine equations, and Baker presents some of these. More about that just below. The proofs follow the so-called hypergeometric method used by Thue in his earlier paper of 1908 (also see just below).

Actually Thue in a later paper from 1918 had implicitly (in terms of diophantine equations) obtained results like (\ref{lb14}). His work implies for example an effective
$$\left|\alpha -  {p \over q}\right|  \geq {c \over q^{6.668}}$$
for $\alpha=17^{1/7}=1.4989\ldots$ thanks to $p_?/q_?=3/2=1.5$; a bit more elegant than (\ref{lb14}).

\medskip
But the ultimate in elegance was then achieved by Baker \cite{6w} also in 1964: namely
\begin{equation}\label{lb15}
\left|2^{1/3}  -  {p \over q}\right|  > {10^{-6}\over q^{2.955}}.
\end{equation}
The consequence for diophantine equations is that all solutions $x,y$ in $\bf Z$ of
\begin{equation}\label{de16}
x^3-2y^3=m
\end{equation}
satisfy
\begin{equation}\label{ub16a}
\max\{|x|,|y|\} \leq (300000|m|)^{23}.
\end{equation}
Even though equations like (\ref{de16}) had been around since at least Thue's work, this was the first time that (\ref{de16}) itself could be solved even in principle. It is similar to the situation around (\ref{mord1}) and (\ref{mord2}) above. 

Thus to solve completely for example
\begin{equation}\label{thue}
x^3-2y^3=1621
\end{equation} 
we only have to check all $x,y$ with absolute value at most say $10^{200}$. And with modern computers this is entirely feasible, using a trick found by Baker and Davenport a few years later (see section 5). 

Actually to prove these results one needs infinitely many elements like $p_?/q_?$ in the set $\bf Q$ of all rationals. An old idea, going back essentially to Hermite, for getting many approximations $p/q$ in $\bf Q$ to a number $\xi$ in $\bf R$ is to construct many approximations $A/B$ in ${\bf Q}(x)$ (preferably of so-called Pad\'e type) to a suitable function $f(x)$. Then one specializes $x$.

The choice $f(x)=(1-x)^{1/3}$ will do for (\ref{lb15}). Thue had already used this, and found that hypergeometric functions turn up. Baker succeeded with (\ref{lb15}) by means of a nice new twist. 

By easy linear algebra we see that for any $r \geq 0$ in $\bf Z$ there are polynomials $A_r(x),B_r(x)$, not both zero, such that 
\begin{equation}\label{pad17}
\phi_r(x)=A_r(x)-(1-x)^{1/3}B_r(x)
\end{equation}
has a zero at $x=0$ of order at least $2r+1$ (these are in fact Pad\'e approximations). They are unique if we normalize $B_r(0)=1$ and then 
$$B_r(x)=F\left({1 \over 3}-r,-r,-2r,x\right)=\sum_{j=0}^rb_{rj}x^j$$
for the Gauss function $F$ (now a polynomial). Here the coefficients 
$$b_{rj}=\prod_{k=0}^{j-1}{((1/3)-r+k)(-r+k) \over (1+k)(-2r+k)}$$ 
are all rational, and seem to involve $3^j$ in their denominators.

As $2^{1/3}=1.2599\ldots$ the choice $p_?/q_?=5/4$ is tempting in the discussion above; and this translates into putting $x=3/128$ in (\ref{pad17}), when $(1-x)^{1/3}={5 \over 4} / 2^{1/3}$. As $\phi_r(x)=cx^{2r+1}+\cdots$ is small, it seems plausible that we obtain a sequence of good rational approximations 
\begin{equation}\label{app18}
{5 \over 4}{B_r(3/128) \over A_r(3/128)}~~~~(r=0,1,2,\ldots)
\end{equation}
to $2^{1/3}$. For example with $r=0,1,2,3,4$ we get
$${5 \over 4},~{635 \over 504},~{96389 \over 76504},~{15240955\over12096754},~{26990767415\over21422586294}$$
(the last being accurate to 19 decimal places). The hope is that these will repel all others as in (\ref{tri}). But to check this we have to estimate denominators in (\ref{app18}). For example
$$B_r\left({3 \over 128}\right)=\sum_{j=0}^rb_{rj}\left({3 \over 128}\right)^j;$$
and now we see that the $3^j$ in the denominators of $b_{rj}$ are cancelled out. It is this ``3-adic'' feature that is Baker's new twist (which led later to a minor industry). It would fail for the slightly more tempting $p_?/q_?=63/50=1.26$.

In the same paper he did similar things with $a^{1/3}$ for $a=17,19,20,37$ and $43$. But, at least up to 1988, the method fails for $a=5$ and so does not lead to a complete solution of say $x^3-5y^3=1621$ (see section 4).

\medskip
In all Baker published no less than five papers in 1964, each highly non-trivial (he did even better with six in 1967). We have already discussed four of them. The fifth is \cite{2w}, which follows the broad principles of \cite{5w} and \cite{6w} above. Now he is dealing with numbers $\xi=\log \alpha$ whose irrationality (or transcendence) is by itself somewhat deep. In fact for algebraic $\alpha \neq 0,1$ this was proved by Lindemann in 1882 (see just below in section 4). We indicate just one consequence in the style of (\ref{lb14}): if $b=10^{12}$ then
$$\left|\log \left({b+1 \over b}\right)  -  {p \over q}\right|  \geq {10^{-121506} \over q^{2.1}}.$$
But he is also able to extend to transcendence measures in the style of ({\ref{tm10}), with $\kappa$ arbitrarily near $n$ in the spirit of (\ref{tm10a}). In the proofs the analogue of (\ref{pad17}) is a system of forms
\begin{equation}\label{pad18}
\phi(x)=A^{(0}(x)+A^{(1)}(x)\log(1-x)+\cdots+A^{(n)}(x)(\log(1-x))^n
\end{equation}
found by Mahler.

\medskip
This category of Baker's work is rounded off by \cite{14w} in 1967, which extends \cite{5w} to linear forms in various 
$$(a/b)^{d_0/d},(a/b)^{d_1/d},\ldots,(a/b)^{d_n/d},$$
again close to best possible in the spirit of (\ref{tm10a}). Instead of quoting an approximation result we mention an attractive application to diophantine equations. Namely if $h=10^{11}+1$ and $l=h^5-1$ then all solutions $x,y,z$ in $\bf Z$ of
\begin{equation}\label{de18}
x^5+ly^5+l^2z^5+5lxyz(xz-y^2)=m
\end{equation}
satisfy
$$\max\{|x|,|y|,|z|\} \leq l^{500}m^2.$$
This may not be quite so elegant as (\ref{de16}) and (\ref{ub16a}); but the exponent of $m$ is much smaller, and crucially this seems to be the first such example ever for equations in three variables (where at the time of writing it is unknown if $x^3+y^3+z^3=114$ - see \cite{sum} - has any solutions at all). Now another system like (\ref{pad18}), also due to Mahler, plays a key role.

\bigskip
This concludes the section on diophantine approximation. Baker did more very important work around (\ref{liou6}), but that deserves an extra section, which now follows. In the sequence of papers just described he seems to be exuberantly flexing his muscles on several different methods, (with hindsight) limbering up for his big achievements. In the early stages he was possibly solving every problem Davenport threw at him; but pretty soon his own motivation took over.

\section{Linear forms in logarithms}\label{sec4}

This title hardly existed before Alan Baker. 

\medskip
The classical theorem of Hermite-Lindemann (1873-1882) is equivalent to the fact that if $\alpha$ is a non-zero algebraic number, and $\log\alpha$ is any non-zero choice of its complex logarithm, then $1$ and $\log\alpha$ are linearly independent over the field $\overline{\bf Q}$ of all complex algebraic numbers. This includes the transcendence of $e,\pi$ and so
$$\int_0^1{{\rm d}x \over x^2+1}={\pi \over 4},$$ 
and of course the numbers $\log\alpha$ mentioned in section \ref{sec3}, as well as the numbers $e^\beta$ for algebraic $\beta \neq 0$.

Similarly the classical theorem of A.O. Gelfond and T. Schneider (1929-1934) is equivalent to the fact that if $\alpha_1,\alpha_2$ are non-zero algebraic numbers, and $\log\alpha_1,\log\alpha_2$ are any choices of logarithms which are linearly independent over $\bf Q$, then they are even linearly independent over $\overline {\bf Q}$. This includes the transcendence of $e^\pi$ and most notably the numbers $2^{\sqrt{2}}$ and $\alpha^\beta$ mentioned by Hilbert (Seventh Problem) in his famous 1900 address to the International Congress of Mathematicians in Paris, as well as numbers like $\log 3/\log 2$. 

Before Baker, nothing was known about the independence of $1,\log \alpha_1,\log \alpha_2$; and practically nothing about that of $\log \alpha_1,\log \alpha_2,\log \alpha_3$, despite the assertions of Gelfond and Linnik in 1948 (repeated in Gelfond's 1960 book) of the enormous importance of the latter problem (for which see later). There had been a paper by Gelfond and N.I. Feldman in 1949 about 
$$b_1\log \alpha_1+b_2\log \alpha_2+b_3\log \alpha_3$$
for integers $b_1,b_2,b_3$; but these coefficients were subject to a restrictive condition which essentially reduced the problem to two logarithms.

\medskip
Then Baker in a sequence of four papers 1966-68 (which we shall call the Linear Forms Quartet) went straight for any number of logarithms and proved (among much else) 

\medskip
\noindent
{\bf Theorem (Baker).} {\it If $\alpha_1,\dots,\alpha_n$ are non-zero algebraic numbers, and $\log\alpha_1,\ldots,\log\alpha_n$ are any choices of logarithms which are linearly independent over $\bf Q$, then 
\begin{equation}\label{logs5}
1,\log\alpha_1,\ldots,\log\alpha_n
\end{equation} 
are linearly independent over $\overline {\bf Q}$.} 

\medskip
The reader may easily construct simple examples of transcendental numbers like $e^{\beta_0}\alpha^\beta$ or $\alpha_1^{\beta_1}\alpha_2^{\beta_2}$ not covered by Hermite-Lindemann or Gelfond-Schneider;  less simple examples are 
\begin{equation}\label{ints}
\int_0^1{{\rm d}x \over x^3+1}~=~{\pi\sqrt{3} \over 9}+{\log 2 \over 3}
\end{equation}
quoted by Siegel in his famous 1949 transcendence monograph, or
$$\int_0^1{{\rm d}x \over x^3-3x-1}~=~\beta_1\log\alpha_1+\beta_2\log\alpha_2+\beta_3\log\alpha_3$$
where 
$$\alpha_1=4-\alpha^2,~\alpha_2=2+\alpha,~\alpha_3=-\alpha+\alpha^2,$$
$$\beta_1={-4-\alpha+2\alpha^2 \over 9},~\beta_2={2-\alpha-\alpha^2 \over 9},~\beta_3={2+2\alpha-\alpha^2 \over 9}$$
and $\alpha^3-3\alpha-1=0$ (with say $\alpha>0$ to provide unique real choices of logarithms).

\medskip
We proceed to examine these four papers I, II, III, IV in detail.

\medskip
Linear Forms I is \cite{11w1} from 1966. It is one of the great number-theoretic papers of the twentieth century. 

\medskip
In fact the main result is a ``linear independence measure'' (analogous to transcendence measure in section \ref{sec3}) for just $\log\alpha_1,\ldots,\log\alpha_n$ with a condition slightly stronger than that in the above theorem. We return to such measures later; it is these refinements that are needed for the applications suggested by Gelfond and Linnik.

\medskip
In an exemplary display of modesty, clarity and foresight, Baker writes 

\medskip
``{\it Finally, as regards the proof of the theorem, our method depends on the construction of an auxiliary function of several complex variables which would seem to be the natural generalisation of the function of a single variable used in Gelfond's original work. The subsequent treatment employed by Gelfond, however, is not applicable in the more general context and so it has been necessary to devise a new technique. Nevertheless it will be appreciated that the argument involves many familiar ideas. The method will probably be capable of considerable development for it applies in principle to many other auxiliary functions apart from the one constructed here.}''

\medskip
This is immediately followed by

\medskip
``{\it The author is grateful to Prof. H. Davenport, who read the original draft of this paper, for his helpful criticism.}''

\medskip
\noindent
which makes it clear that Davenport did not suggest the problem. Indeed the problem, despite its great importance, seems not to have been very well known; for example there is no mention of it in Lang's book on transcendental numbers, which appeared in the same year 1966. It is not referred to in the section ``{\it Einige offene Fragestellungen}'' in Schneider's book of 1957; this section contains his celebrated Eight Problems. And in Siegel's book from 1949 there is no mention of logarithms in connexion with (\ref{ints}). Also in Baker's paper there is no serious reference to any paper published after 1948. He sought out the problem himself, despite all previous transcendence breakthroughs coming through French, Russian and German sources; and solved it in spectacular style.

\medskip
Let us give an idea of his ``new technique''.

We start with $n=2$ and Gelfond-Schneider; slightly changing the notation we have to deduce a contradiction from a relation
\begin{equation}\label{rel}
\beta\log\alpha=\log\alpha'
\end{equation} with $\alpha,\alpha',\beta$ algebraic and $\beta$ irrational. Gelfond constructs an auxiliary function $\Phi(z)$ which is a non-zero polynomial of large degree in 
\begin{equation}\label{fun}
e^{z},e^{\beta z}.
\end{equation}
Note that from (\ref{rel}) the functions (\ref{fun}) take algebraic values at all points
\begin{equation}\label{point}
z~=~s\log\alpha,~~~~~s=0,1,2,\ldots
\end{equation}
and this is true even of their derivatives, thanks to the differential equation for the exponential function. That enables $\Phi$ to be constructed, with algebraic coefficients not all zero, such that
\begin{equation}\label{eqgel}
{{\rm d}^{t} \over {\rm d}z^{t}}\Phi(s\log\alpha)~=~0
\end{equation}
for all non-negative integers $s,t$ in some large range
\begin{equation}\label{ran}
s \leq S,~~t \leq T.
\end{equation}
Thus $\Phi$ has many zeroes (with multiplicity) inside a large disc in $\bf C$ (say centred at the origin).

Now well-known analytic techniques such as the Schwarz Lemma (or Maximum Modulus Principle - see (\ref{schw}) below) show that $\Phi$ must be very small on the whole of this disc, and even on a slightly larger one; and this holds for the derivatives, and even slightly more of them. Then well-known arithmetic techniques show that (\ref{eqgel}) holds for say
\begin{equation}\label{ran+}
s \leq 2S,~~t \leq 2T.
\end{equation}

This step could be then iterated, even indefinitely, to get infinitely many zeroes. In fact a single zero of infinite multiplicity, say at $z=0$, suffices to prove that $\Phi$ must be identically zero, and this leads easily to the required contradiction.

\medskip
How did Baker adapt Gelfond's proof? For simplicity we take $n=3$ and ignore the extra $1$ in (\ref{logs5}), so that now we have to deduce a contradiction from a relation
\begin{equation}\label{rel6}
\beta_1\log\alpha_1+\beta_2\log\alpha_2~=~\log\alpha_3
\end{equation}
with $\beta_1,\beta_2$ algebraic (with a suitable irrationality condition). The auxiliary function $\Phi(z_1,z_2)$ is now a polynomial of large degree in
\begin{equation}\label{fun7}
e^{z_1},e^{z_2},e^{\beta_1z_1+\beta_2z_2}
\end{equation}
which indeed generalize in a fairly natural way (although no-one had previously written them down) Gelfond's $\Phi(z)$ and (\ref{fun}). This time (\ref{rel6}) shows that the functions (\ref{fun7}) take algebraic values at all points
\begin{equation}\label{zero}
(z_1,z_2)~=~(s\log\alpha_1,s\log\alpha_2),~~~s=0,1,2,\ldots
\end{equation}
and this is true also of their partial derivatives. That enables $\Phi$ to be constructed, with algebraic coefficients not all zero, such that
\begin{equation}\label{eqbak}
{\partial^{t_1} \over \partial z_1^{t_1}}{\partial^{t_2} \over \partial z_2^{t_2}}\Phi(s\log\alpha_1,s\log\alpha_2)~=~0
\end{equation}
for all non-negative integers $s,t_1,t_2$ in some large range
\begin{equation}\label{ran8}
s \leq S,~~t_1+t_2 \leq T.
\end{equation}
Thus $\Phi$ has many zeroes (still with a natural concept of multiplicity) inside a large disc now in ${\bf C}^2$.

In 1941 Schneider had used several complex variables for transcendence purposes, and in ${\bf C}^2$ was able to go further because the set of zeroes was part of a structure like ${\bf Z}^2$. (Later Baker expressed it ``{\it ... this type of argument  requires that the points in question form a cartesian product, a condition that can apparently be satisfied only with respect to particular multiply-periodic functions.}'') But in (\ref{zero}) the structure is only $\bf Z$; and indeed the points all lie on a line. This is too thin a set for the usual type of Schwarz Lemma. Indeed the function $(z_1\log\alpha_2-z_2\log\alpha_1)^T$ has the same zeroes as $\Phi$ but it need not be small even on a disc of radius 1 in ${\bf C}^2$.

And to this day, no-one knows how to increase $T$ in (\ref{ran8}) to $2T$ as in (\ref{ran+}). 

Here comes Baker's decisive innovation. It does seem fairly natural to restrict to the line. But Baker considered also the partial derivatives
$$\Phi_{t_1,t_2}=\Phi_{t_1,t_2}(z)~=~{\partial^{t_1} \over \partial z_1^{t_1}}{\partial^{t_2} \over \partial z_2^{t_2}}\Phi(z\log\alpha_1,z\log\alpha_2)$$
however keeping only $t_1+t_2 \leq {1 \over 2}T$.

Clearly ${{\rm d}^{t} \Phi_{t_1,t_2}/ {\rm d}z^{t}}$ can be expressed as a linear combination of the $\Phi_{\tau_1,\tau_2}$ with
$$\tau_1+\tau_2 \leq {1 \over 2}T+t.$$
So as long as we keep also $t \leq {1 \over 2}T$, then $\tau_1+\tau_2 \leq T$ and we see from (\ref{eqbak}) that 
$${{\rm d}^{t} \over {\rm d}z^{t}}\Phi_{t_1,t_2}(s)~=~0$$
for 
$$s \leq S,~~t \leq {1 \over 2}T.$$
These are very similar to (\ref{eqgel}) and (\ref{ran}), firmly back in $\bf C$.

Now the Schwarz Lemma shows that the $\Phi_{t_1,t_2}$ are very small on a slightly large disc; and so the arithmetic gives us (\ref{eqbak}) for
\begin{equation}\label{ran8+}
s \leq 2S,~~t_1+t_2 \leq {1 \over 2}T.
\end{equation}

We can iterate but not indefinitely in any profitable way, due to the loss of multiplicity. 

About this Martin Huxley has commented (1967):

\medskip
	``{\it Where GEL'FOND had one exponential,
	
	ALAN BAKER foresaw the potential.
	
	A strange iteration
	
	Gives zero-inflation.
	
	The right place to start is essential.}''
	
\medskip	
Already this limited iteration was a new sort of difficulty, which Baker overcame by getting just as many zeroes as are needed for the contradiction (along the principle that a polynomial of degree $D$ cannot have $D+1$ zeroes). Here the non-vanishing of a certain Vandermonde determinant suffices.

\medskip
Of course all this extends to any number of variables and
\begin{equation}\label{eqbakn}
{\partial^{t_1} \over \partial z_1^{t_1}}\cdots{\partial^{t_n} \over \partial z_n^{t_n}}\Phi(s\log\alpha_1,\ldots,s\log\alpha_n)~=~0.
\end{equation}

\medskip
So Baker's main achievement was to introduce several complex variables, reduce them to a single variable along a line, introduce the magic trick with multiplicities, and supply the missing estimates for zeroes.

\medskip
Actually he was helped in this last ``zero estimate'' by Gelfond himself. In a 1935 paper he had obtained a linear independence measure for two logarithms. This amounts to assuming that $\beta\log\alpha-\log\alpha'$ is very small rather than zero as in (\ref{rel}). That adds some technical complications but more significantly it means that here too one cannot iterate (\ref{ran+}) indefinitely; hence the need for some sort of zero estimate. 

To state Gelfond's result we have at last to introduce heights. Recall that the degree of an algebraic number $\alpha$ is the smallest $d$ in any equation (\ref{alg5}) satisfied by $x=\alpha$. We can then assume the coefficients have no common factor, and we define the height by
$${\mathcal H}(\alpha)=\max\{|b_0|,|b_1|,\ldots,|b_d|\}$$
(we are using $\mathcal H$ to distinguish it from a more modern version that need not concern us here - see however the end of this section).

Thus Gelfond was able to show that for any $\kappa>5$ (and later for any $\kappa>3$) there is effective $c>0$ depending only on $\alpha_1,\alpha_2$ and the degrees of $\beta_1,\beta_2$ such that
\begin{equation}\label{lbgel}
|\beta_1\log \alpha_1+\beta_2\log\alpha_2| \geq ce^{-(\log B)^\kappa},
\end{equation}
where $B=\max\{2,{\mathcal H}(\beta_1),{\mathcal H}(\beta_2)\}$, provided $\log\alpha_1,\log\alpha_2$ are linearly independent over $\bf Q$ and their principal values are chosen.

Accordingly in Linear Forms I Baker could prove an analogous
\begin{equation}\label{lbbak}
|\beta_1\log \alpha_1+\cdots+\beta_n\log\alpha_n| \geq ce^{-(\log B)^\kappa}
\end{equation}
for any $\kappa>n+1$, where now $B=\max\{2,{\mathcal H}(\beta_1),\ldots,{\mathcal H}(\beta_n)\}$, except that now $\log\alpha_1,\ldots,\log\alpha_n,2\pi i$ have to be linearly independent over $\bf Q$.

\medskip
In the second paper Linear Forms II \cite{11w2} from 1967 Baker obtained (\ref{lbbak}) for any $\kappa>2n+1$ when only $\log\alpha_1,\ldots,\log\alpha_n$ are linearly independent over $\bf Q$. This was crucial to an application to an effective improvement of Liouville's (\ref{liou6}) (see later).

\medskip
In Linear Forms III \cite{11w3}, also from 1967, Baker improved the above result of \cite{11w2} for any $\kappa>n$. He also made a start on lower bounds for the ``inhomogeneous form''
$$\beta_0+\beta_1\log \alpha_1+\cdots+\beta_n\log\alpha_n.$$
This finally yields the theorem quoted above.

\medskip
In the last paper of the Quartet, Linear Forms IV \cite{18w} of 1968, Baker starts by noting that the various positive constants $c$ in his lower bounds, although of course effective, are very small indeed, which could lead to practical computing difficulties in the applications. He gives the following version leading to explicit $c$ which are not so small. Suppose for some $\delta$ with $0<\delta \leq 1$ that there exist $b_1,\ldots,b_n$ in $\bf Z$ with absolute values at most $H$ such that
\begin{equation}\label{ineq}
0<|b_1\log \alpha_1+\cdots+b_n\log\alpha_n| <e^{-\delta H}.
\end{equation}
Then
\begin{equation}\label{hb}
H<(4^{n^2}\delta^{-1}d^{2n}\log A)^{(2n+1)^2}
\end{equation}
where $d \geq 4$ is an upper bound for the degrees of $\alpha_1,\ldots,\alpha_n$ and $A \geq 4$ an upper bound for their heights.

Implicit in this paper is the following associated result (see as well \cite{37w} p.3 and p.8), also the first of its kind. If $\log \alpha_1,\ldots,\log\alpha_n$ are this time linearly dependent over $\bf Q$ then there is a ``relatively bounded'' relation
\begin{equation}\label{rb1}
0=b_1\log \alpha_1+\cdots+b_n\log\alpha_n
\end{equation}
for $b_1,\ldots,b_n$ in $\bf Z$ with
\begin{equation}\label{rb2}
0<\max\{|b_1|,\ldots,|b_n|\} <(4^{n^2}d^{2n}\log A)^{(2n+1)^2}.
\end{equation}
The idea here was greatly exploited in later work on isogenies between elliptic curves and between abelian varieties and consequences of that (see (\ref{rb3}) below).

\medskip
Shortly afterwards Baker wrote two important papers \cite{15w} and \cite{16w} which improve further the Quartet results; but since the main emphasis there is on diophantine equations we shall postpone our description to the corresponding category. It is also \cite{15w} that contains the first effective improvement on Liouville's (\ref{liou6}).

\medskip
Similar remarks hold for his 1971 paper \cite{26w}, published in the top journal Annals of Mathematics, in connexion with class numbers. 

\medskip
We turn now to the later developments of (\ref{lbbak}). These may seem rather technical, but they are specifically aimed at marvellous applications, to be described in sections 5 and 7.

\medskip
We start with his paper \cite{61w} with Harold Stark, also in the Annals and 1971.

\medskip
Already in \cite{15w} Baker had observed that the term $(\log A)^{(2n+1)^2}$ in (\ref{hb}) could be replaced by $C(\log A_n)^\kappa$ for any $\kappa>n+1$, where $A_n\geq 4$ is an upper bound for ${\mathcal H}(\alpha_n)$ alone, and $C$ (effective) is now allowed to depend on $\delta,\kappa,d$ as well as $\alpha_1,\ldots,\alpha_{n-1}$. In \cite{61w} they replace (\ref{ineq}) by
\begin{equation}\label{ineq2}
0<|\beta_1\log \alpha_1+\cdots+\beta_n\log\alpha_n| <e^{-\delta H}
\end{equation}
where now $\beta_1,\ldots,\beta_n$ are algebraic numbers also of degree at most $d$. The analogue of (\ref{hb}) is
\begin{equation}\label{bh}
H<\max\{e^{\sqrt{\log B}},C'(\log A_n)^\kappa\}
\end{equation}
now for any $\kappa>1$, where $B\geq 4$ is an upper bound for the heights of $\beta_1,\ldots,\beta_n$ and $C'$ has the same nature as $C$ above.

This indeed technical result could be applied to class number problems (see section 7); but also to general effective improvements of Liouville's (\ref{liou6}), also described later on. But actually much interest attaches to the proof. Already in \cite{23w} (see section 6) Baker had the nice idea of taking $s$ in (\ref{point}) and (\ref{zero}) as rational not just integral, but now they combine this with another important idea, which for simplicity we illustrate with (\ref{point}). Namely the functions (\ref{fun}) take at $z=(1/p)\log\alpha$ the values 
$$\alpha^{1/p},\alpha'^{1/p}$$
and the vanishing of $\Phi$ there implies an algebraic relation between these. Such things can be ruled out by Kummer theory (based on Galois theory) if $p$ is for example a sufficiently large prime (thus $1,2^{1/p},4^{1/p}$ are linearly independent over $\bf Q$ for $p>2$); and this gives another approach to the zero estimates.

\medskip
It is not difficult to guess which author wrote the sentence

\medskip
``{\it However, the latter work involves an appeal to Galois theory which we prefer to avoid.}''

\medskip
\noindent
in a footnote.

\medskip
We now come to the ``Sharpening Trio'' I,II,III (changed at the last moment from ``strengthening'').

\medskip
In Sharpening I \cite{28w} (in memory of Davenport and Sierpi\'nski) from 1972 Baker is able to improve (\ref{bh}) for $\delta=1$ to
$$H<C''\log B\log A_n$$
provided $\beta_1,\ldots,\beta_n$ are in $\bf Z$, where $C''$ is allowed to depend only on $d,\alpha_1,\ldots,\alpha_{n-1}$. Now the main motivation has become effective improvements of Liouville, and a little competition opens up between Baker and Feldman. Indeed the much sharper dependence on $B$ (now essentially best possible) comes from a 1968 paper \cite{fel} of Feldman (who realized how to avoid certain factorials - Thue himself had to do this as well). Also Baker concedes in a footnote that Feldman has recently obtained those effective improvements (for which having $\kappa=1$ in (\ref{bh}) is essential), but immediately adds that he can do it himself with a ``slight generalization'' of this \cite{28w}.

\medskip
This slight generalization is Sharpening II \cite{32w}, from 1973 (dedicated to Siegel). It contains a partial improvement of (\ref{hb}) when $b_n=-1$ in (\ref{ineq}); namely
$$H < C'''\log A_n,$$
where the ever-effective $C'''$ is allowed to depend only on $d,\delta,\alpha_1,\ldots,\alpha_{n-1}$. Indeed this does do the trick for Liouville, and it is customary to declare the competition a draw. Thus we have 

\medskip
\noindent
{\bf Theorem (Baker-Feldman).} {\it Given any algebraic number $\alpha$ of degree $d \geq 3$ there is effective $\kappa<d$ and effective $c>0$ such that
\begin{equation}\label{bf}
\left|\alpha - {p \over q}\right| \geq {c \over q^\kappa}
\end{equation}
for all $p$ and $q \geq 1$ in $\bf Z$.}

\medskip
Regarding this $\kappa$, we pause to mention the 1988 work \cite{63w} of Baker and Cam Stewart. They gave explicit values of $\kappa$ and $c$ for every irrational $\alpha=a^{1/3}$ with $a \geq 2$ in $\bf Z$. We already mentioned a problem with $a=5$ in connexion with (\ref{lb15}); and this they solve with
$$\left|5^{1/3}  -  {p \over q}\right|  ~>~  {10^{-12900}\over q^{2.9999999999998}}$$
(which would have left Estermann's Scotsman without words). The proofs use cubic units (here $41+24\alpha+14\alpha^2$) and an {\it ad hoc} linear form in three logarithms.

\medskip
The Trio was interrupted with \cite{29w} in 1973. Here Baker introduces the now-familiar notation
\begin{equation}\label{lam}
\Lambda=\beta_0+\beta_1\log \alpha_1+\cdots+\beta_n\log\alpha_n
\end{equation}
(back to inhomogeneous) and shows that if $\Lambda \neq 0$ then
\begin{equation}\label{lamlb}
|\Lambda| > (B\log A_n)^{-C''''\log A_n}
\end{equation}
with $A_n$ as above, and now $B\geq 4$ is an upper bound for the heights of $\beta_0,\beta_1,\ldots,\beta_n$ and $C''''$ is allowed to depend effectively on $d,\alpha_1,\ldots,\alpha_{n-1}$. This starts to look more like the modern versions. He gives the application
$$\left|e^\pi - {p \over q}\right| \geq {1 \over q^{C\log\log q}}$$
for all $p$ and $q \geq 4$ in $\bf Z$, with $C$ absolute (that is, depending on no additional parameters); this has never been improved. 

\medskip
The Trio ends with Sharpening III \cite{33w} (in memory of Linnik and Mordell) from 1975. There is now a mini-competition with Stark. Baker shows that when $\beta_0=0$ and $\beta_1,\ldots,\beta_n$ are in $\bf Z$ in (\ref{lam}) still with $\Lambda \neq 0$, then
\begin{equation}\label{lamlb2}
|\Lambda| > B^{-C\Omega\log \Omega}.
\end{equation}
Here $\Omega=\log A_1\cdots\log A_n$ where each $A_i \geq 4$ is an upper bound for the height ${\mathcal H}(\alpha_i)$. Things seem to be getting more complicated; but now the effective $C$ depends only on $n$ and $d$, so we are heading back in the direction of totally explicit estimates like (\ref{hb}). And the temporary neglect of $\alpha_1,\ldots,\alpha_{n-1}$ has ended (this had been done by Stark two years earlier).

Baker gives no applications; but a footnote ``added in proof'' mentions that the result has ``recently'' been used by Robert Tijdeman to deliver what is perhaps the most spectacular application of linear forms in logarithms to diophantine equations; more of that in the next section.

\medskip
With \cite{37w} from 1977 we are finally back to totally explicit estimates. After a detailed account of the history (up to then), Baker shows firstly that for $\Lambda \neq 0$ as in (\ref{lam}) we have
\begin{equation}\label{lamlb3}
|\Lambda| > (B\Omega)^{-C_0\Omega\log \Omega'}
\end{equation}
with $\Omega'=\Omega/\log A_n$ and
$$C_0=(16nD)^{200n};$$
now 
$$D=[{\bf Q}(\alpha_1,\ldots,\alpha_n,\beta_0,\beta_1,\ldots,\beta_n):{\bf Q}]$$
measures degrees in a different way. And secondly that when $\beta_0=0$ and $\beta_1,\ldots,\beta_n$ are in $\bf Z$ the bracketed factor $\Omega$ in (\ref{lamlb3}) can be eliminated to yield
\begin{equation}\label{lamlb4}
|\Lambda| > B^{-C_0\Omega\log \Omega'}.
\end{equation}

\medskip
The proof supplements his own formidable techniques with ideas of Tijdeman and Alf van der Poorten (about the new $\Omega'$) and Tarlok Shorey (about $C_0$).

\medskip
And there Baker's extraordinary solo work on linear forms in logarithms ended. We have already mentioned a paper with Stewart in 1988 but it was not until 1993 that he did any more work on (\ref{lamlb4}).

\medskip
This was \cite{64w} with Gisbert W\"ustholz. The main result was expressed with a more modern version of height and also took into account non-principal values of the logarithms, but we give only the consequence for principal values. Again we are in the so-called ``rational case'' when $\beta_0=0$ and $\beta_1,\ldots,\beta_n$ are in $\bf Z$. Then they improve (\ref{lamlb4}) to
\begin{equation}\label{lamlb5}
|\Lambda| > B^{-C_1\Omega}
\end{equation}
with 
$$C_1=(16nD)^{2n+4}.$$
The elimination of $\log\Omega'$ is due mainly to the zero estimates of W\"ustholz arising from his more general work on group varieties.

\medskip
Incredibly sharp as all these results are, we should note that there is much room for conjectural improvement; for example one believes that the product of logarithms $\Omega$ can be replaced by the sum (see section 8).

Much more recently Vesselin Dimitrov and Philipp Habegger \cite{DH} have made a conjecture which would have valuable consequences in the theory of dynamical systems of algebraic origin. In our present framework it runs as follows. Given $\epsilon>0$ and $T \geq 1$, there exists $C=C(\epsilon,T)$ such that 
$$|b_1\log \alpha+b_2\log(-1)| > \exp(-CD(DB)^\epsilon)$$
for all algebraic $\alpha$, not a root of unity, with $|\alpha|=1$ and ${\mathcal H}_M(\alpha) \leq T$. Here the notation is by now familiar except for ${\mathcal H}_M$, which is the standard modern version of height (or Mahler measure) defined using the conjugates of $\alpha$. The crux is the very good dependence on $D$.

\bigskip
The theory of Linear Forms in Logarithms (or logarithmic forms for short) continues to be one of the two main tools for solving diophantine problems in number theory (or indeed several other branches of mathematics), with a key emphasis on effectivity. The other is the Subspace Theorem due to Schmidt with important extensions by Hans Peter Schlickewei, which is partly non-effective.

\section{Diophantine equations}\label{sec5}
We already had several examples of these, such as (\ref{thue}), and we even mentioned the year 1621; but the so-called Pell equation
\begin{equation}\label{pell}
x^2-410286423278424y^2 = 1
\end{equation}
was encountered in the solution of the ``Cattle Problem'' of Archimedes (see \cite{var} for example).

\medskip
In general it is a question of finding a solution with the variables in $\bf Z$ (which any idiot can do for Pell above) or more ambitiously finding all solutions (there are infinitely many for Pell above, all but two of them safely out of the reach of idiots).

\medskip
To warm up let us show that there are at most finitely many positive integers $r,s$ with 
\begin{equation}\label{sunit}
3^r-2^s=1621
\end{equation} 
(the equation $3^r-2^s=1$, with a musicological background, was already solved by Levi ben Gershon in 1343 - see \cite{sim} for example). Assume $H=\max\{r,s\}$ is large, say $H \geq H_0$. Then $3^r$ is relatively close to $2^s$ and of course
\begin{equation}\label{59a}
3^r2^{-s}=1+{1621 \over 2^s}.
\end{equation}
Taking logarithms shows that $r\log 3-s\log 2$ is very small, and we get without trouble (\ref{ineq}) with $n=2,\alpha_1=3,\alpha_2=2,b_1=r,b_2=-s,\delta=1/2$ for an easily computable $H_0$. So (\ref{hb}) with $d=4,A=4$ gives $H < 10^{132}$ and ben's your uncle.

\medskip
But for most of this section we shall consider only polynomial equations. It was Thue in 1909 who first proved that equations like $x^3-2y^3=1621$ in (\ref{thue}) have at most finitely many solutions, as a consequence of his improvement on Liouville's (\ref{liou6}). But we already observed the ineffectivity in that his method did not allow all the solutions to be found. We also saw Baker overcome this problem with (\ref{ub16a}) for $x^3-2y^3=m$ in (\ref{de16}). However we noted that his method did not work for $x^3-5y^3=m$.

\medskip
It was in \cite{15w} (part I of a Duo published in the Philosophical Transactions of this Society) from 1968 that Baker was able to treat the general equation considered by Thue. Namely, if $f$ is a homogeneous  polynomial in $x,y$ of degree $d \geq 3$ (so Pell is out), irreducible over $\bf Q$ and with coefficients in $\bf Z$, and $m \neq 0$ is in $\bf Z$, then all solutions of
$$f(x,y)=m$$
satisfy
\begin{equation}\label{thueub}
\max\{|x|,|y|\} < Ce^{(\log |m|)^\kappa}
\end{equation}
for any $\kappa>d+1$, where $C$ depends only on $f$ and $\kappa$.

He deduces easily the very first effective improvement on (\ref{liou6}); namely  
\begin{equation}\label{liou66}
\left|\alpha - {p \over q}\right| \geq {ce^{(\log q)^{1/\kappa}}\over q^d},
\end{equation}
so he is already on his way to (\ref{bf}).

The proof uses a new estimate for linear forms in logarithms (what else?), which we already commented on. We sketch the method, because it shows that arbitrarily many logarithms may be needed.

\medskip
For simplicity we take the equation
\begin{equation}\label{thued}
x^d-2y^d=1
\end{equation}
with $d \geq 3$ an arbitrary prime. Pick any $\alpha$ with $\alpha^d=2$. Then factorizing the left-hand side of (\ref{thued}) shows that $x-\alpha y$ is a unit in the field ${\bf Q}(\alpha,\zeta)$ with $\zeta=e^{2\pi i/d}$. Dirichlet's Theorem shows the unit group to have rank $t=d(d-1)/2-1$, generated by say $\eta_1,\ldots,\eta_t$ modulo torsion. Thus we may write
$$x-\alpha y=\theta\eta_1^{a_1}\cdots\eta_t^{a_t}$$
with integer exponents and a root of unity $\theta$.

We have $d$ choices for $\alpha$ at our disposal; but it suffices to take just two more, say $\beta,\gamma$. We deduce analogously
$$x-\beta y=\phi\eta_1^{b_1}\cdots\eta_t^{b_t},~~~x-\gamma y=\psi\eta_1^{c_1}\cdots\eta_t^{c_t}.$$

Now ``Siegel's Identity''
$$(\beta-\gamma)(x-\alpha y)+(\gamma-\alpha)(x-\beta y)+(\alpha-\beta)(x-\gamma y)=0$$
leads to
\begin{equation}\label{monom}
\delta\eta_1^{b_1-c_1}\cdots\eta_t^{b_t-c_t}=1+\epsilon
\end{equation}
(a bit as in (\ref{59a}) above) with 
$$\delta=-{(\gamma-\alpha)\phi \over (\alpha-\beta)\psi},~~~\epsilon={(\beta-\gamma)(x-\alpha y) \over (\alpha-\beta)(x-\gamma y)}.$$

Choosing $\alpha$ as the real $2^{1/d}$ it may be seen that if $x,y$ are very large then $x-\alpha y$ and so $\epsilon$ is very small. Thus taking logarithms in (\ref{monom}) gives an integer $a$ such that
$$\log \delta+(b_1-c_1)\log \eta_1+\cdots+(b_t-c_t)\log \eta_t+2\pi ia$$
is also very small (note that $2\pi i$ is a logarithm of 1). This turns out to be more than enough to contradict (\ref{lbbak}).

A similar treatment of $x^d-2y^d=m$ leads to (\ref{liou66}) with $\alpha=2^{1/d}$.

\medskip
In \cite{16w} (part II) also 1968 he solves the Mordell equation $y^2=x^3+k$ in (\ref{mord1}) as we already described in section \ref{sec2}. He also makes (\ref{thueub}) completely explicit with \begin{equation}\label{thueub+}
C=\exp(d^{\nu^2}H^{\nu d^3})
\end{equation}
where $H$ is the maximum of the absolute values of the coefficients of $f$ and $\nu=32\kappa^2d/(\kappa-d-1)$.

The proofs rely on an explicit version of the linear forms estimate in \cite{15w}. The Mordell equation leads to a Thue equation by the classical method of reduction theory of binary cubic forms (which he also has to make explicit along the way, itself of no small interest), and then the general Thue equation needs some simple estimates for units (which were immediately improved by Siegel himself).

\medskip
And thick and fast they came at last. Mathematically next comes \cite{17w} in 1967 (dedicated to Mordell on his 80th birthday). There Baker treats the ``elliptic'' equation
\begin{equation}\label{ell}
y^2=ax^3+bx^2+cx+d
\end{equation}
for $a,b,c,d$ in $\bf Z$ such that the cubic has distinct zeroes. He proves
$$\max\{|x|,|y|\} < \exp\{(10^6H)^{1000000}\}$$
where $H=\max\{|a|,|b|,|c|,|d|\}$. This looks a bit like (\ref{mord2}) but cannot be deduced from it. Instead Baker uses now explicit versions of reduction theory of binary quartic forms and then applies (\ref{thueub+}).

\medskip
Then in \cite{19w} from 1969 he generalizes to the ``hyperelliptic'' equation
\begin{equation}\label{hyp}
y^2=a_0x^n+\cdots+a_n
\end{equation}
with $n \geq 3$ and the polynomial on the right having at least three simple zeroes, and the ``superelliptic'' equation
\begin{equation}\label{sup}
y^m=a_0x^n+\cdots+a_n
\end{equation}
with $m \geq 3, n \geq 3$ and the polynomial on the right having at least two simple zeroes; where the coefficients $a_0,\ldots,a_n$ are in $\bf Z$ again with $H$ as the maximum of their absolute values.

For (\ref{sup}) he proves the doubly exponential
$$\max\{|x|,|y|\} < \exp\exp\{(5m)^{10}(n^{10n}H)^{n^2}\}$$
and for (\ref{hyp}) the triply exponential ``finite but only just''
\begin{equation}\label{ubhyp}
\max\{|x|,|y|\} < \exp\exp\exp\{(n^{10n}H)^{n^2}\}
\end{equation}
(at which the lecture audience always laughs). We emphasize though that these were the first ever bounds of their kind.

The proofs rely on the calculations of \cite{15w} and \cite{16w} in suitably generalized form.

\medskip
In \cite{55w} with Coates from 1970 the bounds get even bigger, but the problem is somewhat different. They take now a polynomial $F$ in $x,y$ with coefficients in $\bf Z$ which is irreducible over $\bf C$. Then
$$F(x,y)=0$$
defines a curve $\mathcal C$ in ${\bf C}^2$ which has a certain geometric genus $g \geq 0$ in $\bf Z$. They assume in fact $g=1$. This does not restrict the total degree $n$ of $F$. Then for all points $(x,y)$ in ${\bf Z}^2$ on $\mathcal C$ they prove
$$\max\{|x|,|y|\} < \exp\exp\exp\{(2H)^{10^{n^{10}}}\}$$
(audience goes crazy). The proof reduces to (\ref{ell}) (but now with algebraic coefficients) using Riemann-Roch in the classical way. Of course this step must be made effective and explicit; and this is done by Coates in a separate paper of great independent value. They can then follow the arguments of \cite{19w} just above.

They remark of their result:

\medskip
``{\it It does not seem to extend easily, however, to curves of genus $> 1$, and an effective proof in the general case remains an important quest.}''

\medskip
And today no-one has the slightest idea how to extend this even to genus $g=2$, even though the finiteness of integral points is known since Siegel in 1929. Of course the standard forms are (\ref{hyp}) for $n=5,6$ so already covered. But there are many other models; for example with Umberto Zannier we found the simple family
$$x^4-y^3-axy=0$$ 
which does not seem to be effectively solvable.

\medskip
We said back in section 2 that bounds like the above can be efficiently dealt with. This was done for the first time in 1969 by Baker and Davenport \cite{57w}.

\medskip
The numbers $1,3,8,120$ have the property that the product of any two, increased by 1, is a square. Van Lint asked if this can be extended to $1,3,8,120,N$ for some integer $N>120$. Plainly then the property holds just for $1,3,8,N$; and he was able to show that if $N$ exists, then
\begin{equation}\label{van}
N \geq 10^{1700000}
\end{equation}
which makes its existence unlikely. In \cite{57w} they prove that indeed there is no $N$.

Again Martin Huxley (around 1969):

\medskip
	``{\it Each product is one from a square.
	
	How many such numbers be there?
	
	The known set of four
	
	Cannot take any more,
	
	As DAVENPORT, BAKER declare.}''

\medskip
It seems worthwhile to give some details, as they also illustrate the raw power of linear forms in logarithms.

\medskip
Clearly if $N$ exists it has the form $x^2-1$ for $x$ in $\bf Z$; and now the other conditions amount to
\begin{equation}\label{pell2}
y^2-3x^2=-2,~~~~z^2-8x^2=-7
\end{equation}
for $y,z$ in $\bf Z$. These imply $(yz)^2=24x^4-37x^2+14$, of the form (\ref{hyp}); but the bound (\ref{ubhyp}) really may be too big.

Instead we use the Pell theory behind things like (\ref{pell}) to see that the first equation of (\ref{pell2}) forces
\begin{equation}\label{eqy}
y+x\sqrt{3}=(1+\sqrt{3})(2+\sqrt{3})^r
\end{equation}
for some $r$ in $\bf Z$; and the second
\begin{equation}\label{eqz}
z+x\sqrt{8}=(1+\sqrt{8})(3+\sqrt{8})^s
\end{equation}
for some $s$ in $\bf Z$ (there is actually a similar second possibility here which we shall ignore).

It is easy to see that (\ref{eqy}) is very near $2x\sqrt{3}$ and (\ref{eqz}) is very near $2x\sqrt{8}$. Eliminating $x$ we see that $\sqrt{8}(1+\sqrt{3})(2+\sqrt{3})^r$ is very near $\sqrt{3}(1+\sqrt{8})(3+\sqrt{8})^s$. This is also like (\ref{sunit}). Taking logarithms in the usual way we see that
$$b_1\log \alpha_1+b_2\log\alpha_2+b_3\log \alpha_3$$
is very small for
$$\alpha_1=2+\sqrt{3},~~\alpha_2=3+\sqrt{8},~~\alpha_3={\sqrt{3}(1+\sqrt{8}) \over \sqrt{8}(1+\sqrt{3})}$$
$$b_1=r,~~b_2=-s,~~b_3=-1.$$
In fact we have (\ref{ineq}), and then (\ref{hb}) leads easily to
\begin{equation}\label{van2}
r<10^{487}.
\end{equation}
But as $N=x^2-1$ grows exponentially in $r$ by (\ref{eqy}), this far from contradicts (\ref{van}).

More explicitly, (\ref{ineq}) leads to (a stronger form of) say
\begin{equation}\label{ineq2}
|r\theta-s+\phi| < 2^{-r}
\end{equation}
for certain real irrational $\theta,\phi$ whose exact form we can now forget.

\medskip
Here comes the new idea. Suppose in place of (\ref{van2}) we have obtained
\begin{equation}\label{van3}
r<R
\end{equation}
for some huge integer $R$. We pick a largish integer $K$ and then by Dirichlet's (\ref{dir}) we find $q$ with
$$\|q\theta\| \leq {1 \over KR},~~~q \leq KR.$$
It follows easily from (\ref{ineq2}) and (\ref{van3}) that
\begin{equation}\label{mod2}
\|q\phi\| \leq {q \over 2^r}+{1 \over K}.
\end{equation}

Now $q$ was chosen specifically to make $q\theta$ very near an integer, and there seems no reason to think that $q\phi$ is also near an integer (shades of Littlewood). For example we have $\|q\phi\|\geq1/200$ with ``probability'' $99/100$. If we could have picked $K$ such that $\|q\phi\| \geq 2/K$ then by (\ref{mod2}) we get
$$2^r \leq qK \leq K^2R.$$
For fixed $K$ this gives a new bound for $r$ that is logarithmic in $R$, and so it is plausible that (\ref{van3}) will be improved.

And indeed in the stronger form of (\ref{ineq2}) they pick $K$ around $10^{33}$, find $q$ efficiently from the continued fraction of $\theta$, and end up with $r<500$ (thanks to the Atlas Computer Laboratory - and presumably not to the pocket calculator that Baker bought after seeing that Eva McLean had one). This does contradict (\ref{van}) and the thing is done.

\medskip
This section would be incomplete without a mention of the Catalan equation
\begin{equation}\label{cat}
x^r-y^s=1\end{equation}
this time to be solved for $x,y,r,s$ in $\bf Z$ all bigger than 1. Thus it is a combination of (\ref{thue}) and (\ref{sunit}). It is looking for all gaps 1 when we throw in fourth powers, fifth powers etc with the squbes of section \ref{sec2}.

In 1842 Catalan conjectured that the only solution is $3^2-2^3=1$. Progress was very slow until a result of Cassels, who also earlier in 1953 had made the weaker conjecture that there are only finitely many solutions. Despite all the work of Baker, it came as a sensation when Tijdeman in 1976 proved the conjecture of Cassels with an ingenious application of linear forms in logarithms. He used Sharpening III (\ref{lamlb2}), together with his own sharpening (now the knives are really out) of Sharpening I. That was something like (\ref{lamlb2}) with a slightly better dependence on $\log A_n$ and a worse dependence on the other factors of $\Omega$ (but not as good as (\ref{lamlb3}) which appeared soon after).

\medskip
At the time it was considered as a pinnacle of ``Bakery'' (van der Poorten) applied to diophantine equations. Since then Preda ${\rm Mih\breve ailescu}$ proved the original Catalan conjecture in 2004. Shortly afterwards in 2006 he found another proof not using linear forms in logarithms (and no computers either). See for example \cite{bbm}.

\section{Elliptic functions}\label{sec6}
\medskip
This section is short but sweet. 

\medskip
We have seen that the Hermite-Lindemann Theorem essentially concerns the exponential function ${\mathcal E}=e^z$, and the simple differential equation ${\mathcal E}'={\mathcal E}$ plays a role in the proofs. It was Siegel in 1921 who first investigated what happens for a Weierstrass elliptic function $\wp=\wp(z)$. This too satisfies a differential equation, namely $\wp'^2=4\wp^3-g_2\wp-g_3$, and since we are working with algebraic numbers it is customary to assume that the ``invariants'' $g_2,g_3$ are algebraic. Minor technical problems are caused by the non-linearity; but a much more serious difficulty comes from analytic growth. The Schwarz Lemma in its simplest form says that if $f$ is an entire function with $f(0)=0$, then 
\begin{equation}\label{schw}
|f(z_0)| \leq {|z_0| \over R}\sup_{|z|=R}|f(z)|
\end{equation}
for any $z_0$ with $|z_0| \leq R$. However $\wp$ is far from entire. We can make it entire by multiplying by the square of the Weierstrass sigma function $\sigma(z)$. Unfortunately $\sup_{|z|=R}|\sigma(z)|$ grows exponentially in $R^2$, while $\sup_{|z|=R}|e^z|=e^R$; and this restricts the iteration procedure in the proof (there is a similar problem with ``arithmetic growth'').

Nevertheless Schneider in 1937 could prove the elliptic analogues of Hermite-Lindemann and Gelfond-Schneider. We formulated those in terms of logarithms 
\begin{equation}\label{log}
\log\alpha=\int_{1}^\alpha {{\rm d}x \over x},
\end{equation}
and the elliptic analogue is
\begin{equation}\label{elog}
\int_{\infty}^\alpha {{\rm d}x \over \sqrt{4x^3-g_2x-g_3}}
\end{equation}
for which there is no standard ``principal value'' (and not even a standard notation). 

Thus it was natural that Baker should start thinking about elliptic analogues of his Theorem quoted in section \ref{sec4}. If we take $\alpha=1$ in (\ref{log}), then there is a loop of integration which gives $2\pi i$, a ``fundamental period'' of $e^z$. Similarly if we take $\alpha=\infty$ in (\ref{elog}), then there are two loops which give ``fundamental periods'' $\omega_1,\omega_2$ of $\wp(z)$, so that
\begin{equation}\label{per}
\wp(z+\omega_1)=\wp(z),~~~~~\wp(z+\omega_2)=\wp(z).
\end{equation}

\medskip
In \cite{23w} from 1970, also a major work with far-reaching consequences later on, he succeeded in proving that
if $\beta_1,\beta_2$ are algebraic, then the linear form in elliptic logarithms $\beta_1\omega_1+\beta_2\omega_2$ is either zero or transcendental.

The analogues of the functions (\ref{fun7}) are
\begin{equation}\label{efun7}
\wp(z_1),\wp(z_2),\beta_1z_1+\beta_2z_2
\end{equation}
and the analogues of the zeroes (\ref{zero}) are
\begin{equation}\label{ezero}
(z_1,z_2)~=~\left(\left(s+{1 \over 2}\right)\omega_1,\left(s+{1 \over 2}\right)\omega_2\right),~~~s=0,1,2,\ldots
\end{equation}
(the $1/2$ just to avoid poles).

However the growth problems mentioned above mean that these zeroes are not sufficiently numerous to do the iterative steps. As hinted in section \ref{sec4}, Baker was forced to use more 
$$s+{p \over q}$$
(for $p=1,\ldots,q-1$) in (\ref{ezero}); such things had not occurred in any earlier transcendence proofs. Then not only must the range of $s$ be increased step by step, as in (\ref{ran8+}); but simultaneously also the range of $q$. 

This brought new problems; for example $\wp(\omega_1/q)$, although still an algebraic number, will have large degree growing like $q^2$ as $q$ increases. Nevertheless Baker was able to get enough zeroes for a contradiction, although the Vandermonde argument had to be ratcheted up a bit.

In fact he proved something more general, with $\omega_1$ as a period of a Weierstrass $\wp_1$ and $\omega_2$ a period of a second Weierstrass $\wp_2$. Schneider had already determined when $\beta_1\omega_1+\beta_2\omega_2$ can be zero.

\medskip
Then in \cite{20w} published in 1969 he considered also not just periods of $\wp$ but also the corresponding quasi-periods $\eta_1,\eta_2$ of the associated Weierstrass zeta function $\zeta(z)$; these satisfy
$$\zeta(z+\omega_1)=\zeta(z)+\eta_1,~~~~~\zeta(z+\omega_2)=\zeta(z)+\eta_2$$
like (\ref{per}). He showed that if $\beta_1,\beta_2,\gamma_1,\gamma_2$ are algebraic then $\beta_1\omega_1+\beta_2\omega_2+\gamma_1\eta_1+\gamma_2\eta_2$ is either zero or transcendental. Already Schneider in 1937 had considered the numbers $\beta_1\omega_1+\gamma_1\eta_1$ and $\beta_2\omega_2+\gamma_2\eta_2$ separately. Now Baker had to deal with new algebraic numbers like $\zeta(\omega_1/q)-\eta_1/q$. He could treat $\wp_1,\wp_2$ also here. In the paper he points out that this implies the transcendence of the sum of the circumferences of two ellipses with algebraic ``axes-lengths''.

\medskip
Finally the note \cite{22w} from 1970 gives a lower bound for $|\beta_0+\beta_1\omega_1+\beta_2\omega_2|$ like (\ref{lbbak}) which has the amusing consequence
$$|\wp(n)| \leq Ce^{(\log n)^\kappa}$$
for any $n \geq 1$ in $\bf Z$.

\section{Class numbers}\label{sec7}
Now we come to more pinnacles of Bakery. 

\medskip
Inside the field $\bf Q$ we have the ring $\bf Z$, and inside $\bf Z$ sit the prime numbers, which it is now convenient to take as
$$\pm2,\pm3,\pm5,\pm7,\pm11,\pm13,\pm17,\pm19,\pm23,\pm29,\ldots,\pm163,\ldots,\pm1979339339,\ldots;$$
and in the present context one calls them ``irreducibles''. Euclid more or less proved that there is unique factorization in the sense that

\medskip
(1) Every $n \neq 0,\pm1$ in $\bf Z$ is a product of irreducibles.

(2) If two such products $p_1\cdots p_k,q_1\cdots q_l$ are equal then $k=l$.

(3) And in (2), possibly after a permutation, we have $p_i=\pm q_i$ for all $i$.

\medskip
Inside a number field $K$ we have similarly a natural ring $\mathcal O$ together with irreducibles, and one can hope for the analogues of (1),(2),(3), except that the ``$\pm$'' must be interpreted in terms of units as in the discussion around (\ref{thued}). And indeed if this were true for all $K$ then ``Fermat's Last Theorem'' would be relatively easy to prove. But it fails already for ${\bf Q}(\sqrt{-5})$.

The nearest approach to these analogues is through ideal classes. One works with ideals $\mathcal A$, which are just (non-zero) ${\mathcal O}$-modules in $\mathcal O$. Two ideals $\mathcal {A,A'}$ are said to be equivalent if there are non-zero $\alpha,\alpha'$ in $\mathcal O$ with $\alpha'{\mathcal A}=\alpha\mathcal A'$. The set of equivalence classes is finite, and its cardinality is the ``class number'' $h$ of $K$.

It is elementary that $h=1$ is the same as the analogues of (1),(2),(3).

When $K$ is quadratic over $\bf Q$ and imaginary, all this was known to Gauss, who in his famous {\it Disquisitiones Arithmeticae} of 1801 calculated many $K$ with small given class number (he actually worked with quadratic forms instead). He conjectured that his lists were complete, but that rigorous proofs seem to be very difficult: ``{\it Demonstrationes autem rigorosae harum observationum perdifficiles esse videntur.}''. Now $K={\bf Q}(\sqrt{-d})$ for some unique positive squarefree integer $d$, and we write $h(d)$ for the class number. His list of $d$ for $h(d)=1$ was
$$1,2,3,7,11,19,43,67,163$$
so nine fields $K$. The last $d$ here corresponds to Euler's discovery that the expression $x^2+x+41$ takes prime values for $x=0,1,2,\ldots,39$, and also to the fact that
$$e^{\pi\sqrt{163}}=262537412640768743.999999999999250072597198185688\ldots$$
is so close to an integer.

In 1934 Heilbronn and Linfoot proved that there is at most one more value of $d$; and if it exists then 
\begin{equation}\label{hl}
d>\exp(10^7).
\end{equation}
But it was not until Baker's work in 1966 that this could be in principle ruled out. In his great paper \cite{11w1} he remarked that Gelfond and Linnik in 1948 had reduced the problem using Kronecker's Limit Formula to a weaker form of his (\ref{lbbak}) with $n=3$ logarithms. So that did it. On 22nd December 1966 J-P. Serre \cite{cor} wrote to John Tate: 

\medskip
``{\it Et Davenport dit qu'un jeune Anglais a prouv\'e l'inexistence du 10-i\`eme corps quadratique imaginaire principal. Pas mal d'un coup!}''

\medskip
\noindent
to which Tate replied on 13th January 1967 in apparent contradiction:

\medskip
``{\it The guy who proved the inexistence of the 10th imaginary quadratic field with $h=1$ is ${\rm American}$, not English. A student of Lehmer, whose thesis proved something like $|D| \geq 10^{10^7}$ if $D<-163$ and $h_D=1$. Now he has done better! He talked recently at MIT about it. I won't take time to try to give the idea now, although I think I could if I tried. Stark (that's his name) says he sees no way to generalize the thing, in order to determine all $D$ with $h(D)=2$, for example. Too bad.}''

\medskip
\noindent
Then Serre on 24th June 1967:

\medskip

``{\it Je reviens \`a Baker: les anglais sont tr\`es excit\'es par son travail, notamment par le fait qu'il donne des bornes calculables pour les probl\`emes d'approximation diophantienne qu'il consid\`ere.}''

\medskip
Actually in \cite{27w} from 1971 he improved the work of Gelfond and Linnik by showing that there is effective $C$ such that if $d>C$ exists (he later calculated $C=10^{20}$), then $21d$ is also squarefree and 
\begin{equation}\label{h1}
\left|h(21d)\log\left({5+\sqrt{21} \over 2}\right)-{32 \over 21}\pi\sqrt{d}\right|<e^{-\pi\sqrt{d}/100}.
\end{equation}
This involves just two logarithms (including $2\pi i$), and classical estimates for $h(21d)$ then show that the inequality can be dealt with by Gelfond's much older bound (\ref{lbgel}). A strange turn of events! Later Baker produced a second inequality like (\ref{h1}) and then eliminated $\pi\sqrt{d}$ to obtain (\ref{ineq}) for $n=2$ for which (\ref{hb}) leads to $d<10^{500}$. So by (\ref{hl}) ``the tenth man'' $d$ does not exist (except as a novel by Graham Greene set in Paris not Vienna).

\medskip
As hinted above, events turn even stranger when it is realized that Stark \cite{sta0} in 1967 had given an independent proof that the tenth man does not exist; that was based on work of Heegner from 1952 which was, at least in some circles, also accepted as a proof -- and then more widely discounted until it was finally recognized as essentially correct.

\medskip
Tate refers to difficulties with class number $h=2$, which for general $K$ is the same as saying that the analogue of (1) or (3) fails but (2) survives (Carlitz). And indeed doing $h(d)=2$ was much tougher. Here Gauss's list was
$$5,6,10,13,15,22,35,37,51,58,91,115,123,187,235,267,403,427.$$

We have already referred to Baker's paper \cite{26w} from 1971; and indeed here is where he showed how to settle this in principle.

\medskip
Already in \cite{21w} from 1969 he had used two logarithms to prove that $d < 10^{500}$ if it is not 3 modulo 8. So  he could assume $d$ is 3 modulo 8. In that case Gauss's theory of genera shows that $d=pq$ for primes $p$ congruent to 1 modulo 4 and $q$ congruent to 3 modulo 4. If $q \leq d^{1/4}$ then also \cite{21w} suffices. So he could also assume that $q>d^{1/4}$.

He then proves an analogue of (\ref{h1}). Namely there is effective $C$ such that if $d>C$ exists, then $21d$ and $21q$ are squarefree and
\begin{equation}\label{h2}
\left|h(21d)\log\left({5+\sqrt{21} \over 2}\right)+h'h(21q)\log \eta-{64 \over 21}\pi\sqrt{d}\right|<e^{-\pi\sqrt{d}/10},
\end{equation}
where $h'$ is the class number of the real quadratic field ${\bf Q}(\sqrt{21p})$ and $\eta$ is a fundamental unit there. An essential difference with (\ref{h1}) is that one of the logarithms, namely $\log \eta$, is not fixed. So nothing from the Quartet can be used, except possibly IV. But there the coefficients were rational; and anyway the exponent $(2n+1)^2=49$ in (\ref{hb}) is far too large. The troublesome $\eta$ might have ${\mathcal H}(\eta)$ as large as $p^{\sqrt{p}}$ so $\log {\mathcal H}(\eta)$ as large as $d^{3/8}$. As $H$ is around $d^{1/2}$ we would need an exponent $\kappa< 4/3$. Baker was able to get any $\kappa>1$, also for algebraic coefficients $\beta_1,\beta_2,\beta_3$ (compare (\ref{bh}) above), thereby settling $h(d)=2$. This was another motivation for the Sharpening Trio. It may be noted that Stark \cite{sta} independently settled the problem, and their papers could not have appeared more simultaneously, with Baker's occupying pages 139 to 152 and Stark's 153 to 173 in volume 94 of the Annals.

\medskip
To this day no-one has been able to treat $h(d)=3$ by similar methods. However the subsequent work of Goldfeld, Gross and Zagier (see \cite{gold} for example) enables in principle any class number to be treated effectively.

\medskip
We already referred to the note \cite{21w}. This was used in a paper \cite{59w} with Andrzej Schinzel (in memory of Davenport), also from 1971, where they make a contribution to the (still unsolved) problem of showing that Euler's list
$$1, 2, 3, 4, 5, 6, 7, 8, 9, 10, 12, 13, 15, 16, 18, 21, 22, 24,\ldots, 408, 462, 520, 760, 840, 1320, 1365, 1848$$
of 65 ``numeri idonei'' is complete. These were defined in terms of binary quadratic forms but they are also the positive integers not of the form $xy+yz+zx$ for distinct positive integers $x,y,z$.

\section{Abcology}\label{sec8}
This refers to the $abc$ Conjecture, formulated in 1985 by Joseph Oesterl\'e and the writer. Namely, for any $\kappa>1$ there is $C=C(\kappa)$ such that for any non-zero $a,b,c$ in $\bf Z$ with no common factor and
$$a+b+c=0$$
we have
\begin{equation}\label{abc1}
\max\{|a|,|b|,|c|\} \leq CS^\kappa
\end{equation}
for the ``squarefree kernel''
$$S = \prod_{p | abc}p \geq 2,$$
the product being taken over all primes $p>0$.

\medskip
In \cite{49w} from 1998 Baker proposes a version
\begin{equation}\label{abc2}
\max\{|a|,|b|,|c|\} \leq C_0S{(\log S)^s \over s!}
\end{equation}
where $s = \sum_{p | abc}1$ and now $C_0$ is absolute. It is easily seen that $s \leq L=[C_0'\log S/\log\log S]$ (if $S \geq 3$) for $C_0'>0$ absolute, and then Stirling gives (for $L \leq \log S$)
$${(\log S)^s \over s!} \leq {(\log S)^{L} \over L!} \leq \left({e\log S \over L}\right)^{L};$$
so (\ref{abc2}) strengthens (\ref{abc1}) - but not so much as to contradict known lower bounds.

He suggested an even stronger version and showed that it is equivalent to a certain lower bound for
$$\min\{1,|\Lambda|\}\prod_p\min\{1,p|\Lambda|_p\}$$
for the linear form in logarithms
$$\Lambda=b_1\log a_1+ \cdots+b_n\log a_n$$
with $b_1,\ldots,b_n,a_1 \geq 1,\ldots,a_n \geq 1$ in $\bf Z$. Here $|\Lambda|_p$ is defined in a natural $p$-adic way. He remarks that this lower bound suggests that the expression $\Omega=\log A_1\cdots\log A_n$ in (\ref{lamlb5}) could be replaced by $\Gamma=\log A_1+\cdots+\log A_n$. He even speculates that $|\Lambda| \geq e^{-C_1\Gamma}B^{-C_1\omega}$, where
$$\omega=\sum_{p | a_1\cdots a_n}1 \leq C_1'{\Gamma \over \log \Gamma}$$
for $C_1,C_1'$ absolute.

\medskip
But despite the extraordinarily precise works of Yu Kunrui on $p$-adic linear forms in logarithms it seems that this approach to $abc$ remains for future generations to explore.

\medskip
Later on in \cite{50w} from 2004 Baker stuck his neck out by conjecturing (on extensive computer evidence) that (\ref{abc2}) holds for $C_0=6/5$.

\medskip
See also his paper \cite{51w} from 2007 for related remarks.

\section{Miscellaneous}\label{sec9}
Here we comment on five papers that do not so neatly fall into the above categories.

\medskip
In \cite{8w} from 1965 he considers algebraic functions $f=f(x)$ satisfying $P(x,f)=0$, where $P(x,y)$ has the special form
$$P(x,y)=1+\left(1+{p \over q}x\right)y+xQ(x,y)$$
for $Q(x,y)$ in ${\bf Z}[x,y]$ and coprime $p,q$ in $\bf Z$ with $q \geq 2$. As $P(0,-1)=0$ but $(\partial P/\partial y)(0,-1) \neq 0$, imposing $f(0)=-1$ gives a unique power series
$$f(x)=-1+a_1x+a_2x^2+a_3x^3+\cdots.$$
He shows that $q^na_n$ is in $\bf Z$ and also prime to $q$. This is in accordance with the well-known Eisenstein Theorem about denominators; but the point is that there is no $q_?\geq 1$ in $\bf Z$ with $q_?^{n-1}a_n$ in $\bf Z$ for all $n \geq 1$, so that the Theorem is best possible in this respect, however we choose $Q$. The simplest choice $Q=0$ gives of course
$$f(x)=-1+{p \over q}x-{p^2 \over q^2}x^2+{p^3 \over q^3}x^3+\cdots.$$

\medskip
In \cite{9w} from 1966 he extends the Pad\'e systems used in \cite{5w}, \cite{6w} (see (\ref{pad17}) above), \cite{2w} (see (\ref{pad18}) above) and (later) \cite{14w}. Namely let $\omega_0,\ldots,\omega_k$ be in $\bf C$ distinct modulo $\bf Z$ and let $\rho_{rs} ~ (r=0,\ldots,k; s=0,\ldots,l)$ be in $\bf Z$ with
$$1 \leq \rho_{rl} \leq \cdots \leq \rho_{r0}~~~(r=0,\ldots,k).$$
Then if $P_{rs}(x)$ is in ${\bf C}[x]$ with degree at most $\rho_{rs}-1$, the function
$$\sum_{r=0}^k\sum_{s=0}^lP_{rs}(x)(1-x)^{\omega_r}(\log(1-x))^s$$
cannot have a zero of order at least $\sum_{r=0}^k\sum_{s=0}^l\rho_{rs}$ at $x=0$ unless all $P_{rs}(x)=0$. He mentions (without proof) that such things can be used to prove the transcendence of $e^\pi$ (see section 4 above).

\medskip
In \cite{12w} from 1967 he extends a famous result of P\'olya to functions of several variables. As the former played an important part in the history of transcendence, and the latter a key role in linear forms in logarithms, we should say a few words.

\medskip
P\'olya shows that if $f(z)$ is an entire function such that $f(0),f(1),f(2),\ldots$ are all in $\bf Z$, then a relatively slow growth rate
$$|f(z)| \leq Ce^{\theta |z|}$$
for some $\theta < \log 2$ and some $C$ forces it to be a polynomial. This led about 20 years later to the Gelfond-Schneider Theorem (compare (\ref{fun}) and (\ref{point}) above). The example $f(z)=2^z$ shows it would be false for $\theta=\log 2$.

Baker proves a similarly sharp result for $f(z_1,\ldots,z_n)$ with values in $\bf Z$ at $z_1,\ldots,z_n=0,1,2,\ldots$. Here it suffices that
$$|f(z_1,\ldots,z_n)| \leq Ce^{\theta (|z_1|+\cdots+|z_n|)},$$
again for some $\theta < \log 2$, to force a polynomial. 

He refers to Schneider's 1941 paper mentioned in section 4. He also gives generalizations involving partial derivatives
$${\partial^{t_1} \over \partial z_1^{t_1}}\cdots{\partial^{t_n} \over \partial z_n^{t_n}}f(z_1,\ldots,z_n)$$
as in (\ref{eqbakn}), so this paper seems inextricably linked with Linear Forms I (which appeared just one year earlier).

\medskip
The paper \cite{53w} from 1973, with Bryan Birch and Eduard Wirsing, uses Linear Forms II and III to prove several results about a problem of Chowla, for example that the values of
$$L(1,\chi)=\sum_{n=1}^\infty {\chi(n) \over n},$$
as $\chi$ runs over all non-principal characters modulo $q$, are linearly independent over $\bf Q$ provided $q$ is prime to its image $\phi(q)$ under the Euler function.

\medskip
The paper \cite{56w} from 1975 with Coates is diophantine approximation of a different sort. Let $p,q$ be coprime in $\bf Z$ with $p>q \geq 2$. It is trivial that
$$\left\|\left({p \over q}\right)^k\right\| \geq {1 \over q^k}$$
for all $k \geq 1$ in $\bf Z$. Mahler showed that for any $\epsilon>0$ there is $c=c(p,q,\epsilon)>0$ such that
$$\left\|\left({p \over q}\right)^k\right\| \geq {c \over q^{\epsilon k}}.$$
For example there are at most finitely many $k$ with
$$\left\|\left({3 \over 2}\right)^k\right\| <{3^k+2^k \over 4^k-2^k};$$
this inequality is significant as (for $k \geq 5$) it is implied by the number $g(k)$ in Waring's Problem NOT being given by
\begin{equation}\label{WP}
g(k)=2^k+\left[\left({3 \over 2}\right)^k\right]-2.
\end{equation}
However the $c$ above is ``Thue-ineffective'', so these exceptional $k$ cannot be determined (it is presumed they do not exist).

\medskip
Baker and Coates show that there are effective $\eta=\eta(p,q)<1$ and $c=c(p,q)>0$ such that
$$\left\|\left({p \over q}\right)^k\right\| \geq {c \over q^{\eta k}}$$
for all $k$. However they admit that getting
$$\eta(3,2)<{\log(4/3) \over \log 2}=.415037\ldots,$$
as would be required to find all exceptions to (\ref{WP}), may need fundamentally new ideas.

The proofs use a refinement of the methods of Linear Forms I, Sharpening I,II and Coates's earlier work on $p$-adic linear forms in logarithms (but not with $p$ as above).
\section{Books}\label{sec10}
Before Baker there were three famous books on transcendental numbers written by leaders in the field: Siegel's ``{\it Transcendental numbers}'' (1949), Schneider's ``{\it Einf\"uhrung in die transzendenten Zahlen}'' (1957), and Gelfond's ``{\it Transcendental and algebraic numbers}'' (1960). The dates reflect a certain irregular growth in the subject corresponding to various breakthroughs. Also Lang's ``{\it Introduction to transcendental numbers''} (1966) was influential. 

\medskip
His own ``{\it Transcendental number theory}'' \cite{43w} from 1975 is a worthy successor to these (although the writer well remembers Baker's wail of agony over a particularly visible printing error on page 85). It formed the substance of an essay for the Adams Prize (1972) of Cambridge University.

After an introductory chapter containing probably some of the world's shortest proofs (one page for the transcendence of $e$, two pages for Lindemann-Weierstrass) he launches straight into linear forms in logarithms, and the subsequent four chapters give an excellent account of this topic and its various applications. They are followed by a chapter on elliptic functions including his own work. The next chapter describes Schmidt's subspace extension of Roth's Theorem. There follows a chapter on Mahler's Classification including Schmidt's beautiful proof that ``{\it $T$-numbers do exist}''. The next two chapters contain some of his own work in \cite{T} and \cite{7w}, and the book finishes with two chapters about the topic of algebraic independence. One is an account of Shidlovsky's fundamental work on $E$-functions and the other contains for example the  proof, due independently to Dale Brownawell and Michel Waldschmidt, that at least one of $e^e$ and $e^{e^2}$ is transcendental. 

Nearly all the chapters contain the first account in book form of major results.

The book was eagerly devoured as soon as it appeared. The concision of the first chapter sets the tone for the rest, which makes it not so easy for beginners. It ran into a second edition in 1990, and at the time of writing a reprint is planned. 

\medskip
Some of the material, now updated, found its way into the first three chapters of his book ``{\it Logarithmic forms and diophantine geometry}'' \cite{66w} with W\"ustholz in 2007. 

The other chapters are then concerned more with developments for general commutative algebraic groups $G$. The simplest examples of these are the additive group ${\bf G}_{\rm a}$ and the multiplicative group ${\bf G}_{\rm m}$. The Quartet results can be naturally formulated in terms of ${\bf G}_{\rm m}^n$ and ${\bf G}_{\rm a}\times{\bf G}_{\rm m}^n$. The next simplest example is an elliptic curve $E$, and Baker's result in \cite{23w} involves in a similar way ${\bf G}_{\rm a}\times E^2$ or even ${\bf G}_{\rm a}\times E_1\times E_2$. And given $E$ there is a unique $G=G(E)$ sitting inside an exact sequence
$$0 \longrightarrow {\bf G}_{\rm a}\longrightarrow G \longrightarrow E \longrightarrow 0$$
which is however not isomorphic to ${\bf G}_{\rm a}\times E$; similarly his paper \cite{20w} involves (implicitly) $G(E_1)\times G(E_2)$ and even a quotient by a line inside.

The fourth chapter supplies the basic theory of such general $G$.

As already mentioned, Baker had to work a bit in \cite{23w} to get enough zeroes for a final contradiction. This indeed was the main obstacle to generalizing even to three elliptic curves. It was overcome by W\"ustholz, and the fifth chapter is devoted to the zero estimate or ``multiplicity estimate'' that did the trick. The resulting extension of Baker's works on linear independence to general $G$ was W\"ustholz's ``Analytic Subgroup Theorem''; and this is the topic of the sixth chapter.

The seventh chapter is concerned with linear independence measures, both improved versions for linear forms in logarithms (as in \cite{64w} above) and versions for general $G$. Some consequences of the latter more in the realm of diophantine geometry (due to W\"ustholz and the writer) are also discussed, such as the existence of ``small''  isogenies between abelian varieties and polarizations of abelian varieties, as well as effective versions of Serre's Open Image Theorem for elliptic curves, the Tate Conjecture for abelian varieties and some of Faltings's Finiteness Theorems.

We pause to make clearer the connexion between isogenies and (\ref{rb1}),(\ref{rb2}). Suppose for example that $E$ and $E'$ are elliptic curves with periods $\omega_1,\omega_2$ and $\omega_1',\omega_2'$ respectively as in (\ref{per}). An isogeny between $E,E'$ leads to relations
\begin{equation}\label{rb3}
\beta \omega_1'=b_{11}\omega_1+b_{12}\omega_2,~~~~\beta \omega_2'=b_{21}\omega_1+b_{22}\omega_2
\end{equation}
as in (\ref{rb1}). The correct analogue of (\ref{rb2}) then leads to an upper bound on the degree of a connecting isogeny.

The book closes with concise accounts of other diophantine topics such as Schmidt's Subspace Theorem and the Andr\'e-Oort Conjecture; the latter takes the reader very near to a great deal of current work.

\medskip
Baker had already published his ``{\it A concise introduction to the theory of numbers}'' \cite{41w} in 1984. This attractive volume is much more elementary and suitable for a first course in Number Theory in general.

\medskip
He then developed \cite{41w} into ``{\it A comprehensive course in number theory}'' \cite{52w} in 2012. This is over twice as long and goes up to the graduate level.

\medskip
Two conferences that Baker organised led to published proceedings: the first (Cambridge 1976) edited by Baker and the writer (and between the two of us we refereed practically all the 16 papers there), mentioned in \cite{37w}; and the second (Durham 1986) edited by Baker alone (where a less intensive programme than usual produced more mathematical advances than usual), mentioned in \cite{63w}.

\section{Fields Medal}\label{sec11}
This medal was for many years the supreme prize in mathematics and still manages very well to hold its own among various lifetime awards. It is given only to people aged at most 40, two or three or four of them every four years at the International Congress. At the time of writing, there have been 60 winners since it began in 1936. 

Baker  won the prize at Nice in 1970, just 13 days after his 31st birthday. See \cite{28ww} for the Proceedings Volume containing his address. The citation was:

\medskip
``{\it Generalized the Gelfond-Schneider theorem (the solution to Hilbert's seventh problem). From this work he generated transcendental numbers not previously identified.''}

\medskip
\noindent
(and indeed at the time of writing that is the only accomplishment listed in his somewhat minimal Wikipedia entry). Paul Tur\'an in his talk at Nice wrote of this: 

\medskip
``{\it The analytic prowess displayed by Baker could hardly receive  a higher testimonial.}'' 

\medskip
\noindent
However any reader of this memoir will see that Baker went far beyond that. And Tur\'an himself recognized the fact in his closing words: 

\medskip
``{\it To conclude, I remark that his work exemplifies two things very convincingly. Firstly, that beside the worthy tendency to start a theory in order to solve a problem it pays also to attack specific difficult problems directly. Particularly is this the case with such problems where rather singular circumstances do not make it probable that a solution would fall out as an easy consequence of a general theory. Secondly, it shows that a direct solution of a deep problem develops itself quite naturally into a healthy theory and gets into early and fruitful contact with other significant problems of mathematics. So, let the two different ways of doing mathematics live in peaceful coexistence for the benefit of our science.}''

\medskip
Baker was presented to Pr\'esident Pompidou at the Elys\'ee Palace and later wrote:  

\medskip
``{\it ... and was quite impressed by his feat of memory when making a welcoming scientific address}''. 

\medskip
He was the third British winner of the Fields Medal (after Klaus Roth 1958 and Michael Atiyah 1966).

\bigskip
Alan Baker single-handedly transformed several areas of number theory. He achieved a major breakthrough in transcendence and applied it to obtain a new and important large class of transcendental numbers (opening the way to the subsequent discovery of several other such classes); developed quantitative versions and applied them to the effective solutions of many classical diophantine equations as well as the first effective improvement on Liouville's 1844 result on diophantine approximation and the resolution of the celebrated Gauss Conjectures of 1801 on class numbers, not only $h=1$ but also $h=2$, of imaginary quadratic fields; and started the study of extensions to elliptic curves (opening the way to later generalizations to abelian varieties and commutative group varieties and in turn their applications to old and new problems in diophantine geometry).

\bigskip
Despite this, his cousin describes him as ``extremely modest''. This is confirmed by another extract from the funeral tribute of Eva McLean:

\medskip
``{\it His most striking characteristic, however, was a genuine modesty. He really believed that his brilliant achievement was merely down to hard work and determination. At the same time he spoke in great awe of those mathematicians whom he regarded as the truly greats.}''
\vglue 1cm
\centerline{ACKNOWLEDGEMENTS}
\bigskip
I thank the Royal Society for allowing me to see their files. I am grateful to Trinity College, particularly to the Librarian Nicolas Bell about Whewell's Court accommodation; he also notes that Alan Baker's personal papers and Fields Medal have been deposited in the College Library.

In addition I have much profited from conversations and correspondence with Sheldon Baker, B\'ela Bollob\'as, William Chen, Yuval Flicker, Martin Huxley, Chris Morley, Sir Michael Pepper, Heather Rechtman, Klaus Schmidt, Cam Stewart, Sir Martin Taylor, Rob Tijdeman, Michel Waldschmidt, Michael Grae Worster, and Gisbert W\"ustholz (also for sending me a complete list of papers) on the personal aspects of Alan Baker's life. But in this respect my greatest debt is of course to Eva McLean.

And on the scientific side I thank too Roger Baker, Enrico Bombieri, Bob Vaughan, and Umberto Zannier, as well as two referees.
\vglue 2cm

\end{document}